\documentclass{amsart}

\usepackage{amsmath}
\usepackage{amsthm}
\usepackage{amsopn}
\usepackage{amssymb}
\usepackage[all]{xy}
\usepackage{graphicx}
\usepackage{lscape}

\newcommand{\abs}[1]{\lvert #1 \rvert}

\newcommand{\br}[1]{\overline{#1}}

\newcommand{\td}[1]{\widetilde{#1}}

\newcommand{\ZZ}{\mathbb{Z}}
\newcommand{\NN}{\mathbb{N}}
\newcommand{\RR}{\mathbb{R}}

\newcommand{\HH}{\mathbb{H}}
\newcommand{\QQ}{\mathbb{Q}}
\newcommand{\WW}{\mathbb{W}}
\newcommand{\FF}{\mathbb{F}}
\newcommand{\GG}{\mathbb{G}}
\newcommand{\MS}{\mathbb{S}}

\newcommand{\calJ}{\mathcal{J}}
\newcommand{\calM}{\mathcal{M}}
\newcommand{\KO}{\mathit{KO}}
\newcommand{\TMF}{\mathit{TMF}}
\newcommand{\tmf}{\mathit{tmf}}

\theoremstyle{definition}
 \newtheorem{thm}[equation]{Theorem}
 \newtheorem{cor}[equation]{Corollary}
 \newtheorem{lem}[equation]{Lemma}
 \newtheorem{prop}[equation]{Proposition}

 \newtheorem{rmk}[equation]{Remark}

 \newtheorem{conj}[equation]{Conjecture}
\newtheorem*{thm*}{Theorem}
\newtheorem*{cor*}{Corollary}
\newtheorem*{lem*}{Lemma}
\newtheorem*{prop*}{Proposition}
\newtheorem*{defn*}{Definition}
\newtheorem*{ex*}{Example}
\newtheorem*{exs*}{Examples}
\newtheorem*{rmk*}{Remark}
\newtheorem*{claim*}{Claim}

\numberwithin{equation}{subsection}
\numberwithin{figure}{subsection}

\DeclareMathOperator{\ord}{ord}
\DeclareMathOperator{\Hom}{Hom}
\DeclareMathOperator{\End}{End}
\DeclareMathOperator{\Res}{Res}

\DeclareMathOperator{\aut}{Aut}
\DeclareMathOperator{\Sta}{Stab}
\DeclareMathOperator{\Sub}{Sub}
\DeclareMathOperator{\Map}{Map}
\DeclareMathOperator{\spec}{spec}
\DeclareMathOperator{\stack}{stack}
\DeclareMathOperator{\Set}{Set}
\DeclareMathOperator{\spf}{spf}
\DeclareMathOperator{\Frob}{Frob}
\DeclareMathOperator*{\holim}{holim}

\DeclareMathOperator*{\colim}{colim}

\title{Buildings, elliptic curves, and the $K(2)$-local sphere}
\author[Mark Behrens]{Mark Behrens$\sp 1$}
\address{
Department of Mathematics \\
Massachusetts Institute of Technology \\
Cambridge, MA 02139, U.S.A.}
\date{\today}

\begin{document}
\begin{abstract}
We investigate a dense subgroup $\Gamma$ of the second Morava stabilizer
group given by a certain group of quasi-isogenies of a supersingular
elliptic curve in characteristic $p$.  
The group $\Gamma$ acts on the Bruhat-Tits building for
$GL_2(\QQ_\ell)$ 
through its action on the $\ell$-adic Tate module.  
This action has
finite stabilizers, giving a small resolution for the
homotopy fixed point spectrum $(E_2^{h\Gamma})^{hGal}$ by spectra of
topological modular forms.
Here, $E_2$ is a version of Morava $E$-theory and 
$Gal = Gal(\br{\FF}_{p}/\FF_p)$.  
\end{abstract}

\maketitle

\tableofcontents

\footnotetext[1]{The author is partially supported by the NSF.}

\section{Introduction}

\subsection{Background}

Fix a prime $p$.  One systematic way of understanding the $p$-local stable 
homotopy groups of a finite 
complex $X$ is to study its chromatic tower, given by the
inverse system
$$ X_{E(0)} \leftarrow X_{E(1)} \leftarrow X_{E(2)} \leftarrow \cdots. $$
Here $X_{E(n)}$ is Bousfield localization with respect to the Johnson-Wilson
spectrum $E(n)$.
The chromatic convergence theorem of Hopkins and Ravenel
\cite{RavenelOrange} states that this tower converges in the sense that
$$ X \simeq \holim_n X_{E(n)} $$
for all $p$-local finite complexes $X$,
and that the induced filtration on the homotopy groups of $X$ is
exhaustive.
The chromatic program for understanding the stable homotopy of $X$ begins
with understanding the filtration quotients of this tower.  This is
equivalent to studying the localizations $X_{K(n)}$ with respect to
Morava $K$-theory.  A very nice
summary of this process may be found in the introduction of \cite{GHMR}.

We confine our attention to the case $X = S$, the sphere spectrum.  
Morava \cite{Morava} developed a method of understanding the 
layers $S_{K(n)}$,
which was strengthened by work of Hopkins-Miller \cite{Rezk},
Goerss-Hopkins \cite{GoerssHopkins}, Devinatz-Hopkins
\cite{DevinatzHopkins}, and Davis \cite{Davis}. 
We briefly summarize their work.
Let $E_n$ be (maximally unramified) Morava $E$-theory.  
It is a complex
orientable spectrum whose associated formal group is the Lubin-Tate
universal deformation of the Honda height $n$ formal group $H_n$ over
$\br{\FF}_p$.  
Our notation is unconventional: $E_n$ is usually taken with respect to the
finite field $\FF_{p^n}$ instead.
We choose to work over $\br{\FF}_p$ because any two
height $n$ formal groups are isomorphic over $\br{\FF}_p$ 
\cite[A2.2.11]{Ravenel}.
Let $\MS_n$ denote the Morava stabilizer group $\aut(H_n)$.
Let $\GG_n$ denote the larger group of automorphisms which are allowed
to act non-trivially on the ground field $\br{\FF}_p$
$$ \GG_n = \aut_{/\FF_p}(H_n) = \MS_n \rtimes Gal. $$
Here $Gal$ is the Galois group $Gal(\br{\FF}_p/\FF_p)$. 
The spectrum $E_n$ is an $E_\infty$-ring spectrum which is 
a continuous $\GG_n$-spectrum, and there is an equivalence
$$ S_{K(n)} \simeq E_n^{h\GG_n} \simeq (E_n^{h\MS_n})^{hGal}. $$

Computationally, it has proven easier to work with homotopy fixed point spectra
$E_n^{hF}$ for \emph{finite} subgroups $F$ of the Morava stabilizer group.
For example, when $n = 1$, the $p$-complete real $K$-theory spectrum 
$\KO_p$ is equivalent to
the homotopy fixed point spectrum $(E_1^{hC_2})^{hGal}$.
Choose $\ell$ to be a topological generator of the group 
$\ZZ_p^\times/\{\pm 1\}$.
The $J$-theory spectrum is given as the fiber
\begin{equation*}
J \rightarrow \KO_p \xrightarrow{\psi^\ell - 1} \KO_p
\end{equation*}
where $\psi^\ell$ is the $\ell$th Adams operation.
Adams-Baird and Ravenel \cite{Bousfield}, \cite{Ravenel84}
proved that there is an equivalence $S_{K(1)} \simeq J$.  Thus the
$K(1)$-local sphere admits a complete description in terms of $K$-theory.

Goerss, Henn, Mahowald, and Rezk \cite{GHMR} gave a similar decomposition
of $S_{K(2)}$ at the prime $3$.  In their work, certain spectra related to
the Hopkins-Miller spectrum of topological modular forms ($\TMF$) played the
role that $K$-theory played in chromatic level $1$.  
This decomposition shed considerable light on the very difficult
$K(2)$-local computations of Shimomura and Wang \cite{ShimomuraWang}.
However, the decomposition was
produced by means of obstruction theory and computation, and the attaching
maps in their decomposition were not identified explicitly.  The
decomposition was also specific to the prime $3$.

\subsection{A higher analog of the $J$-theory spectrum}

It is natural to ask if there is an analog of the $J$-theory spectrum for
chromatic level $2$ which is built out topological modular forms 
in a manner similar to
the way in which the $J$-theory spectrum is built out of $K$-theory.
In \cite{BehrensK(2)}, motivated by \cite{GHMR} and 
\cite{MahowaldRezklevel3}, we introduced a spectrum $Q(\ell)$ 
as the totalization of a 
semi-cosimplicial spectrum
\begin{equation}\label{diag:Q(ell)}
\TMF \Rightarrow \TMF \times \TMF_0(\ell) 
\Rrightarrow \TMF_0(\ell). 
\end{equation}
The spectrum $\TMF$ is the Hopkins-Miller
spectrum of topological modular forms, and the spectrum $\TMF_0(\ell)$ is
an analogous spectrum associated to the congruence subgroup
$\Gamma_0(\ell)$ of $SL_2(\ZZ)$.
Diagram~\ref{diag:Q(ell)} 
has a very natural \emph{abstract} construction in terms of moduli of certain
diagrams of isogenies of elliptic curves.
We refer the reader to
\cite{BehrensK(2)} for this construction, which relies on 
unpublished work of Hopkins, Miller, and their 
collaborators.  However, in \cite{BehrensK(2)} we also gave a
$K(2)$-local
construction for
the case $p = 3$ and $\ell = 2$ which only used 
the Goerss-Hopkins-Miller Theorem.
In Section~\ref{sec:Q(ell)const} we extend this $K(2)$-local construction 
to all
primes $p$ and $\ell$.

The $J$-theory spectrum may also be regarded as the homotopy fixed
point spectrum
$$ J = (E_1^{h\pm\ell^\ZZ})^{hGal} $$
where $\pm\ell^\ZZ$ is the dense subgroup of the Morava stabilizer
group $\MS_1 \cong \ZZ_p^\times$ generated by $\ell$.  The dense subgroup
$\pm\ell^\ZZ$ is the subgroup of the group of 
automorphisms of the multiplicative
formal group $\widehat{\GG}_m$ which is generated by the 
quasi-isogenies of the multiplicative group $\GG_m$ of degree a 
power of $\ell$.

In this paper we aim to give a similar homotopy fixed point construction 
of the spectrum
$Q(\ell)_{K(2)}$.
To this end we define a subgroup $\Gamma$ of $\MS_2$ generated by isogenies
of a supersingular elliptic curve of degree a power of $\ell$.  For
appropriate choices of $\ell$, Tyler Lawson and the author
\cite{BehrensLawson} have shown that 
this subgroup
is dense in $\MS_2$ (it is dense in an index $2$
subgroup if $p=2$).  
There is an extension $\Gamma_{Gal}$ of $\Gamma$ of the form
$$ 1 \rightarrow \Gamma \rightarrow \Gamma_{Gal} \rightarrow \sigma^\ZZ
\rightarrow 1 $$
where $\sigma^\ZZ$ is the dense subgroup of $Gal$ generated by the
Frobenius.  The group $\Gamma_{Gal}$ is dense in $\GG_2$ (respectively, 
an index $2$
subgroup if
$p=2$).
Let $E(\Gamma)$ denote the homotopy
fixed point spectrum
$$ E(\Gamma) = E_2^{h\Gamma_{Gal}}. $$
Because $\Gamma_{Gal}$ is dense in $\GG_2$, one expects that $E(\Gamma)$ is
closely related to the the $K(2)$-local sphere $S_{K(2)} = E_2^{h\GG_2}$.
The precise conjecture is explained in Section~\ref{sec:mainconj}.

We shall discuss the following.
\begin{enumerate}
\item The group $\Gamma$ acts on the Bruhat-Tits building for
$GL_2(\QQ_\ell)$ with finite stabilizers.

\item This action gives a presentation of the group $\Gamma$ 
in terms of
the category of supersingular curves over $\br{\FF}_p$.

\item The action of $\Gamma$ on the building induces a decomposition of 
$E(\Gamma)$ in terms of $K(2)$-local topological modular forms.

\item This decomposition induces an equivalence $Q(\ell)_{K(2)} \simeq
E(\Gamma)$.
\end{enumerate}
The remainder of the introduction is devoted to a more
detailed discussion the results and organization of this paper.

\begin{rmk} 
There is work by other authors which bears some similarity to the contents
of this paper.
\begin{itemize}
\item Gorbounov, Mahowald, and Symonds \cite{GMS}
produced dense amalgamated products of
finite subgroups of the Morava stabilizer group $\MS_{p-1}$.  
Our dense
subgroups appear to differ from theirs in the overlapping case of $p = 3$.

\item In the case of $p=3$, 
the computations of Gorbounov, Siegel, and Symonds \cite{GSS}
reflect algebraically an analog of Conjecture~\ref{conj:mainconj}.

\item Andrew Baker \cite{Baker} has shown that for $p > 3$ the 
$E_2$-term of the 
ANSS for $S_{K(2)}$ can be computed as the cohomology of the 
\emph{$p$-completed} groupoid of supersingular elliptic curves and isogenies.
\end{itemize}
\end{rmk}

\subsection{The subgroup $\Gamma$}
For the remainder of this paper assume we are given a fixed supersingular 
elliptic curve
$C$ over $\br{\FF}_{p}$ and a prime $\ell$ coprime to $p$.  Our work will
turn out to be independent of the choice of supersingular curve --- the 
choice is
tantamount to choosing a basepoint in a connected category
supersingular curves.  For convenience we shall insist that $C$ admits a
definition over $\FF_p$ (for every prime $p$ such a supersingular curve
exists \cite{Waterhouse}).

Since $C$ is supersingular, the
formal completion $C^\wedge$ of $C$ at the identity is isomorphic to the
Honda height $2$ formal group $H_2$.  One may regard $C^\wedge$ as the
$p$-divisible group $C[p^\infty]$.

Our intention is to study the simultaneous action of the endomorphism ring
$\End(C)$ on the $p$-torsion and $\ell$-torsion of $C$.
Let $\Gamma \subset \End(C) \otimes \QQ$ be 
the group of quasi-isogenies of $C$ with degree 
equal to a power of $\ell$.  Then we have the following diagram.
\begin{equation}\label{diag:maindiag}
\xymatrix@C+3.5em{
*+[r]{\underset{\quad}{M_2(\ZZ_\ell)} \cong \End(C[\ell^\infty])} 
\ar@{^{(}->}[d] & 
\: \underset{\quad}{\End(C)} \: 
\ar@{_{(}->}[l]^-{(-)^\wedge_\ell} \ar@{^{(}->}[r]_{(-)^\wedge_p} 
\ar@{^{(}->}[d] &
\underset{\quad}{\End(C[p^\infty])}
\\
*+[r]{M_2(\QQ_\ell)} & 
\; \End(C)[1/\ell] \; \ar@{_{(}->}[l] \ar@{^{(}->}[ru] 
\\
*+[r]{\overset{\quad}{GL_2(\QQ_\ell)}} \ar@{_{(}->}[u] &
\quad \overset{\quad}{\Gamma} \quad 
\ar@{_{(}->}[l] \ar@{_{(}->}[u] \ar@{^{(}->}[r]_{\text{dense}} &
\overset{\quad}{\aut(C[p^\infty]) = \MS_2} \ar@{^{(}->}[uu] 
}
\end{equation}
Tate \cite{Tate},\cite{MilneWaterhouse} 
proved that the top inclusions are actually the 
$\ell$ and $p$-completions of the endomorphism ring, respectively.
In \cite{BehrensLawson}, we proved the following theorem.

\begin{thm}[Behrens-Lawson \cite{BehrensLawson}]\label{thm:Gamma}
Let $\ell$ be a topological generator of $\ZZ_p^\times$ ($\ZZ_2^\times/\{
\pm 1\}$ if $p = 2$).
For $p > 2$, the group $\Gamma$ is dense
in $\MS_2$.  For $p = 2$, the group $\Gamma$ is dense in the index $2$
subgroup $\td{\MS}_2$ which is the kernel of the composite
$$ \MS_2 \xrightarrow{N} \ZZ_2^\times \rightarrow
(\ZZ/8^\times)/\{1,\ell\}. $$
\end{thm}

\subsection{The building}
We shall denote $X^{ss}$ to be the set of isomorphism classes of 
supersingular elliptic curves over $\br{\FF}_p$.
Let $X_0^{ss}(\ell)$ be the set of isomorphism classes of 
pairs $(C',H)$ where $C'$ is a
supersingular elliptic curve and $H$ is a $\Gamma_0(\ell)$-structure (a 
subgroup of order 
$\ell$ contained in $C'(\br{\FF}_p)$).
Given a pair $(C',H) \in X_0^{ss}(\ell)$, let $\aut(C',H)$ be the group of
automorphisms $\phi$ of $C'$ such that $\phi(H) = H$.

The group $\Gamma$ naturally acts on the $\ell$-adic Tate module 
$V_\ell(C)$, giving an inclusion into
the group $GL_2(\QQ_\ell)$.  Let $\calJ'$ be the Bruhat-Tits building for
$GL_2(\QQ_\ell)$. 
$\calJ'$ is a $2$-dimensional contractible simplicial
complex on which $GL_2(\QQ_\ell)$ acts.  
The induced $\Gamma$ action on $\calJ'$ has finite stabilizers, 
which are given naturally by certain groups automorphisms of supersingular
elliptic curves.

\begin{thm*}[Theorem~\ref{thm:calJ'}]
The building $\calJ'$ is $\Gamma$-equivariantly homeomorphic to
the geometric realization of a $\Gamma$-equivariant semi-simplicial complex
$\calJ'_\bullet$.
The simplices of $\calJ'_\bullet$ are given as follows:
\begin{align*}
\calJ'_0 & = \coprod_{C' \in X^{ss}} \Gamma/\aut(C'), \\
\calJ'_1 & = \coprod_{(C',H) \in X_0^{ss}(\ell)} \Gamma/\aut(C',H)
\; \amalg \; \coprod_{C' \in X^{ss}} \Gamma/\aut(C'), \\ 
\calJ'_2 & = \coprod_{(C',H) \in X^{ss}_0(\ell)} \Gamma/\aut(C',H).
\end{align*}
\end{thm*}

\subsection{A resolution of $E(\Gamma)$}

The semi-simplicial complex of
Theorem~\ref{thm:calJ'} gives rise
to the following semi-cosimplicial construction.

\begin{prop*}[Proposition~\ref{prop:Ehgammacosimp}]
There is a semi-cosimplicial $E_\infty$-ring spectrum of the form
$$
\begin{array}{c}
\underset{C' \in X^{ss}}{\prod} E_2^{h\aut(C')}
\end{array}
\Rightarrow 
\begin{array}{c}
\underset{C' \in X^{ss}}{\prod} E_2^{h\aut(C')} \\
\times \\
\underset{(C',H) \in X^{ss}_0(\ell)}{\prod}
E_2^{h\aut(C',H)}
\end{array}
\Rrightarrow 
\begin{array}{c}
\underset{(C',H) \in X^{ss}_0(\ell)}{\prod} 
E_2^{h\aut(C',H)}
\end{array}
$$
which totalizes to give the homotopy fixed point spectrum $E_2^{h\Gamma}$.
\end{prop*}

As we shall explain in Section~\ref{sec:topmodforms},
the $K(2)$-localizations of $\TMF$ and $\TMF_0(\ell)$ are given as 
products of the following homotopy fixed point
spectra:
\begin{align*}
\TMF_{K(2)} & \simeq \left( \prod_{C' \in X^{ss}} E_2^{h\aut(C')}
\right)^{hGal}, \\
\TMF_0(\ell)_{K(2)} & 
\simeq \left( \prod_{(C',H) \in X^{ss}_0(\ell)} 
E_2^{h\aut(C',H)} \right)^{hGal}.
\end{align*}
These constructions, due to Hopkins, Miller, and their collaborators, have
not yet appeared in the literature.  Section~\ref{sec:topmodforms} 
may be regarded
as a self-contained construction of the $K(2)$-local versions of these
spectra of topological modular forms.

It turns out, upon taking Galois homotopy fixed points, 
that the semi-cosimplicial spectrum given by
Proposition~\ref{prop:Ehgammacosimp}  
is the $K(2)$-localization of the 
semi-cosimplicial spectrum of Diagram~(\ref{diag:Q(ell)}) 
defining $Q(\ell)$.  We therefore have the following
theorem.

\begin{thm*}[Theorem~\ref{thm:Q(ell)}]
The spectra $Q(\ell)_{K(2)}$ and $E(\Gamma)$ are naturally equivalent.
\end{thm*}

We also produce a decomposition of $E_2^{h\Gamma^1_{Gal}}$
(Theorem~\ref{thm:hGamma1}), where $\Gamma^1_{Gal}$ is the subgroup of $\Gamma$
consisting of elements of norm $1$.

\subsection{Relation to $K(2)$-local sphere}\label{sec:mainconj}
While the spectrum $J$ is equivalent to the $K(1)$-local sphere, it
appears that the
$K(2)$-local sphere is built out of the spectrum $Q(\ell)_{K(2)}$ and a 
spectrum
dual to $Q(\ell)_{K(2)}$.  More precisely, we make the following conjecture.

\begin{conj}\label{conj:mainconj}
Let $\ell$ be a generator of $\ZZ_p^\times$ 
(respectively $\ZZ_2^\times/ \{\pm 1 \}$ for $p = 2$).
Then if $p$ is odd, the sequence
\begin{equation}\label{eq:mainconj}
D_{K(2)}Q(\ell) \xrightarrow{D\eta} S_{K(2)} \xrightarrow{\eta} 
Q(\ell)_{K(2)}
\end{equation}
is a cofiber sequence.
Here $\eta$ is the $K(2)$-localization of the unit of the ring spectrum
$Q(\ell)$, and $D_{K(2)}$ denotes the Spanier-Whitehead dual in the
$K(2)$-local category.
At the prime $2$, let $\td{S}$ denote the homotopy fixed point spectrum
$E_2^{h\td{\GG}_2}$.  Here, $\td{\GG}_2 = \td{\MS}_2 \rtimes Gal$ 
is an index $2$ subgroup of $\GG_2$, where $\td{\MS}_2$ is the group of
Theorem~\ref{thm:Gamma}.
Then there is a cofiber sequence
$$ D_{\td{S},K(2)}Q(\ell) \xrightarrow{D\eta} \td{S} \xrightarrow{\eta} 
Q(\ell)_{K(2)} $$
where $D_{\td{S},K(2)}$ denotes the Spanier-Whitehead dual in the category of
$K(2)$-local $\td{S}$-modules.
\end{conj}

\begin{rmk}
There is an
equivalence
$$ S_{K(2)} \simeq \td{S}^{hC_2}. $$
The generator of the group $C_2$ lifts to a torsion-free element of the
group $\GG_2$.  Therefore, the $K(2)$-local sphere at the prime $2$
differs mildly from the spectrum $\td{S}$.
\end{rmk}

\begin{rmk}
Conjecture~\ref{conj:mainconj} hypothesizes that the
sequence~(\ref{eq:mainconj})\emph{extends} to a cofiber sequence.  
There are possibly
many different extensions, and the lack of a natural candidate
represents a major
gap in our understanding of $K(2)$-local homotopy theory.
\end{rmk}

The intuition that something like Conjecture~\ref{conj:mainconj} should be
true is due to Mark Mahowald.  
In \cite{BehrensK(2)}, we proved Conjecture~\ref{conj:mainconj} 
in the case $p=3$ and $\ell = 2$.  The author intends to combine
Theorem~\ref{thm:Gamma} with Theorem~\ref{thm:Q(ell)} to
prove a version of Conjecture~\ref{conj:mainconj} for $p > 3$ in a future
paper.

\subsection{Organization of the paper}

In Section~\ref{sec:endC}, we describe the ring of endomorphisms of the
supersingular curve $C$, and describe its action on the formal group of $C$
and the $\ell$-adic Tate module of $C$.  We define the group
$\Gamma$ and show that it is may be viewed as an $\ell$-arithmetic group
associated to a form of $GL_2$.  We also describe an $SL_2$-variant, and
define an associated dense subgroup $\Gamma^1$ of the norm $1$ subgroup 
$\MS l_2 \subset \MS_2$.

In Section~\ref{sec:J'}, we introduce the building $\calJ'$ for
$GL_2(\QQ_\ell)$ from the point of view of $\ZZ_\ell$-lattices in 
$\QQ_\ell^2$.  We then translate this description to one in terms of
subgroups of $C$ using the Weil pairing.  We immediately deduce the
$\Gamma$-equivariant structure of $\calJ'$.

In Section~\ref{sec:J}, we introduce the building (tree) for
$SL_2(\QQ_\ell)$ and run a similar analysis to that of
Section~\ref{sec:J'}, with the group $\Gamma$ replaced by its norm $1$
counterpart $\Gamma^1$.  We deduce amalgamation formulas for $\Gamma^1$
using Bass-Serre theory.

We begin Section~\ref{sec:topmodforms} with a review of the
Goerss-Hopkins-Miller functor, and a technical discussion of the homotopy
fixed point construction of Devinatz and Hopkins.
We then give an exposition of the 
construction of the spectra $\TMF_{K(2)}$ and $\TMF_0(\ell)_{K(2)}$ of Goerss, 
Hopkins, Miller and their collaborators.

In Section~\ref{sec:Qell}, we give a $K(2)$-local construction of the 
spectrum $Q(\ell)$.
We then show
that this spectrum is naturally equivalent to the spectrum $E(\Gamma)$.  We
end by describing a variant where the group $\Gamma$ is replaced by the
norm $1$ subgroup $\Gamma^1$.

\subsection*{Acknowledgments}

The author would like to thank Daniel Davis, Paul Goerss, Hans-Werner Henn, 
Mike Hopkins, Johan de Jong, Tyler Lawson,
Mark Mahowald, Cathy O'Neil, Charles Rezk, and John Rognes.  
This paper
would not have materialized without the generosity with which they shared
their mathematical knowledge and ideas.  

\section{The ring of endomorphisms of $C$}\label{sec:endC}

Let $\End(C)$ be the ring of endomorphisms of $C$ defined over the
algebraic closure $\br{\FF}_p$.  Let $D$ be the ring of quasi-isogenies
$$ D = \End^0(C) = \End(C)\otimes \QQ. $$
Because $C$ is supersingular, $D$ is the quaternion algebra over $\QQ$ 
ramified at $p$ and
$\infty$ \cite{Silverman}.  The subring $\End(C) \subset D$ is a maximal
order.  
For $v$ a valuation, let $D_v = D \otimes_\QQ \QQ_v$ be the completion of $D$
at $v$.
We shall recall in this section how the 
elliptic curve $C$ gives a very explicit description of these local
algebras.

\subsection{The $\ell$-torsion of $C$}

For any prime $\ell$ different from $p$,
there is a non-canonical isomorphism of groups
$$ C[\ell^\infty] \cong \ZZ/\ell^\infty \times \ZZ/\ell^\infty. $$

Let $T_\ell(C)$ be the $\ell$-adic Tate module.  It is the inverse limit of
the inverse system
$$ C[\ell] \xleftarrow{[\ell]} C[\ell^2] \xleftarrow{[\ell]} C[\ell^3]
\xleftarrow{[\ell]} \cdots . $$
The $\ZZ_\ell$-module $T_{\ell}(C)$ is free of rank $2$.
Since every endomorphism of $C$ restricts to an endomorphism of 
$C[\ell^k]$, we see that $T_{\ell}(C)$ is a module over the ring
$\End(C)$. 
We recall the following fundamental theorem of Tate (the case where $A$ is
an elliptic curve, as well as Corollary~\ref{cor:ltorsion}, 
may be deduced from the work of Deuring \cite{Deuring}).

\begin{thm}[Tate \cite{Tate}]\label{thm:Tate}
Let $A$ is an abelian variety over the finite field $\FF_q$.
Then the natural map
$$ \End(A)\otimes \ZZ_\ell \rightarrow 
\End_{\ZZ_\ell[\Frob_q^{rel}]}(T_\ell(A)) $$
is an isomorphism, where $\Frob_q^{rel}$ is the endomorphism induced by the 
$q$th relative Frobenius.
\end{thm}

For supersingular elliptic curves $C'$, some power of the 
relative Frobenius will
lie in the center of $\End(T_\ell(C'))$, so we have the following corollary.

\begin{cor}\label{cor:ltorsion}
Let $C'$ be a supersingular elliptic curve over $\br{\FF}_p$.  
Then the natural map
$$ \End(C')\otimes \ZZ_\ell \rightarrow \End_{\ZZ_\ell}(T_\ell(C')) $$
is an isomorphism.
\end{cor}

\begin{cor}
The algebra $D_\ell$ is split (isomorphic to $M_2(\QQ_\ell)$).
\end{cor}

The Tate module may be equated with the Pontryagin dual of the 
$\ell$-torsion subgroup $C[\ell^\infty]$
as follows.

\begin{prop}\label{prop:Weil}
The Weil pairing induces a Galois equivariant isomorphism
$$ \td{e}: T_\ell(C) \rightarrow 
\Hom(C[\ell^\infty],{\mu}_{\ell^\infty}) = C[\ell^\infty]^* $$
where ${\mu}_{\ell^\infty}$ is the $\ell$-torsion in the
multiplicative group
$\br{\FF}_p^\times$, and the Galois group acts by conjugation on
$\Hom(C[\ell^\infty], \mu_{\ell^\infty})$.
\end{prop}

\begin{proof}
Recall that the Weil pairing is a bilinear Galois equivariant 
non-degenerate pairing
$$ e_{\ell^k}: C[\ell^k] \times C[\ell^k] \rightarrow {\mu}_{\ell^k}. $$
Non-degeneracy implies that the adjoint homomorphism is an isomorphism.
\begin{gather*}
\td{e}_{\ell^k}: C[\ell^k] \xrightarrow{\cong} \Hom(C[\ell^k], 
{\mu}_{\ell^k}) \\
x \mapsto e_{\ell^k}(x,-)
\end{gather*}
One of the properties of the Weil pairing is that the following diagram
commutes \cite[III.8.1]{Silverman}
$$
\xymatrix@R+3em{
C[\ell^{k+1}] \ar[r]^-{\td{e}_{l^{k+1}}} \ar[d]_{[\ell]} 
& \Hom(C[\ell^{k+1}], {\mu}_{\ell^{k+1}}) \ar[d]^{\iota^*}
\\
C[\ell^k] \ar[r]^-{\td{e}_{l^{k}}}  
& \Hom(C[\ell^k], {\mu}_{\ell^k})
}$$
where $\iota: C[\ell^k] \hookrightarrow C[\ell^{k+1}]$ is the inclusion.
The isomorphism $\td{e}$ is the composite
\begin{align*}
T_\ell(C) & = \lim_k C[\ell^k] \\
& \xrightarrow{\cong} \lim_k \Hom(C[\ell^k], {\mu}_{\ell^k}) \\
& \xrightarrow{\cong} \lim_k \Hom(C[\ell^k], {\mu}_{\ell^\infty}) \\
& \xrightarrow{\cong} \Hom(\colim_{k} C[\ell^k], {\mu}_{\ell^\infty}) \\
& \xrightarrow{\cong} \Hom(C[\ell^{\infty}], {\mu}_{\ell^\infty}).
\end{align*}
\end{proof}

The isomorphism $\td{e}$ induces an $\End(C)$-module structure on 
$\Hom(C[\ell^{\infty}], {\mu}_{\ell^\infty})$.  This action is given
explicitly in the following lemma.

\begin{lem}\label{lem:Weilaction}
Let $\alpha$ be an element of 
$C[\ell^\infty]^* = \Hom(C[\ell^{\infty}], 
\mu_{\ell^\infty})$, and let $\phi$ be an endomorphism of $C$.  
Then the action
of $\phi$ on $\alpha$ is given by pre-composition
$$ \phi \cdot \alpha = \alpha \circ \widehat{\phi} $$
where $\widehat{\phi}$ is the dual isogeny.
\end{lem}

\begin{proof}
This is immediate from the following ajointness property of the Weil pairing
\cite[III.8.2]{Silverman}.  For $x$ and $y$ in $C[\ell^k]$, we have
$$ e_{\ell^k}(\phi (x),y) = e_{\ell^k}(x, \widehat{\phi}(y)). $$
\end{proof}

\subsection{The $p$-torsion of $C$}

Because $C$ is supersingular, it has no non-trivial $p$-torsion points.  
The $p$-divisible group $C[p^\infty]$ is 
entirely formal,
meaning that it coincides with the height $2$ formal group $C^\wedge$.  

The endomorphism ring of $C^\wedge$ is the maximal order of the 
$\QQ_p$-division algebra $D_{p,1/2}$ of Hasse invariant $1/2$.
The following theorem is due to Tate.
\begin{thm}[Tate {\cite{MilneWaterhouse}}]\label{thm:ptorsion}
The natural map
$$ \End(C) \otimes \ZZ_p \rightarrow \End(C[p^\infty]) = \End(C^\wedge) $$
is an isomorphism.
\end{thm}

\begin{cor}
The algebra $D_p$ is non-split (isomorphic to $D_{p,1/2}$).
\end{cor}

\begin{rmk}
The fundamental exact sequence of 
class field theory implies that the local invariants of $D$ must add to
zero.  Therefore, the quaternion algebra $D$ must ramify at infinity,
giving an isomorphism
$$ D_\infty \cong \HH. $$
\end{rmk}

\subsection{The reduced norm}

Let $R$ be a ring.
Consider the degree map 
$$ \mathrm{deg}: \End(C) \rightarrow \ZZ. $$
If we choose an additive basis of $\End(C)$, then the degree map is
expressed by a degree $2$ polynomial in $4$ variables.
The degree map 
extends multiplicatively to a reduced norm
$$ N_R = \mathrm{deg}\otimes R: \End(C) \otimes R \rightarrow R. $$
In particular, $N_\QQ$ 
coincides with the reduced norm of the quaternion algebra $D$.

\subsection{The group scheme $G$}

The various groups which appear in Diagram~(\ref{diag:maindiag}) are
conveniently given as the $R$-points of an affine group scheme $G$ for
various $R$.
Define $G$ to be the scheme
whose $R$-points are given by
\begin{align*} G(R)
& = (\End(C) \otimes R)^\times \\
& = \{ x \in \End(C) \otimes R \: : \: N_R(x) \in R^\times \}. 
\end{align*}
This functor is represented by an affine scheme because $\End(C)$ is free
abelian and the reduced norm is given
by a polynomial.

Corollary~\ref{cor:ltorsion} and Theorem~\ref{thm:ptorsion} 
identify the $R$ points of $G$ for various $R$.

\begin{prop}
We have the following values of the functor $G(-)$ where $\ell$ is prime to
$p$:
\begin{align*}
G(\ZZ) & = \aut(C), \\
G(\QQ) & = D^\times, \\
G(\ZZ_\ell) & = \aut(C[\ell^\infty]) = GL(T_\ell(C)) \cong GL_2(\ZZ_\ell), \\
G(\QQ_\ell) & = GL(T_\ell(C)\otimes \QQ) \cong GL_2(\QQ_\ell), \\
G(\ZZ_p) & = \aut(C^\wedge) = \MS_2, \\
G(\QQ_p) & = D_{p,1/2}^\times, \\
G(\RR) & = \HH^\times.
\end{align*}
\end{prop}

\subsection{The group $\Gamma$}\label{sec:Gamma}

We now fix $\ell$ to be a topological 
generator of the group $\ZZ_p^\times$ (respectively 
$\ZZ_2^\times/\{\pm 1\}$).
Define $\End_\ell(C)$ to be the monoid of endomorphisms of $C$ with degree
a power of $\ell$.
Let $\Gamma$ be the group completion of this monoid.

\begin{lem}\label{lem:gpcomplete}
The group $\Gamma$ is given by inverting the element $[\ell]$ of the monoid
$\End_\ell(C)$.  That is to say, the natural map
$$ \End_\ell(C)[\ell^{-1}] \rightarrow \Gamma $$
is an isomorphism.
\end{lem}

\begin{proof}
We simply need to show that $\End_\ell(C)[\ell^{-1}]$ contains
inverses for every $\phi \in \End_\ell(C)$.  If $\phi$ has degree $\ell^k$,
then the dual isogeny $\widehat{\phi}$ has the property
\cite[III.6.2]{Silverman}
$$ \phi \widehat{\phi} = \widehat{\phi} \phi = [\ell^k]. $$
Therefore, the element $\ell^{-k} \cdot \widehat{\phi} \in
\End_\ell(C)[\ell^{-1}]$ is an inverse for $\phi$.
\end{proof}

The group $\Gamma$ is therefore the group of 
quasi-isogenies of 
$C$ with degree a power of $\ell$.  
Alternatively, $\Gamma$
is given by $G$ as
$$ \Gamma = G(\ZZ[1/\ell]). $$
There are inclusions
\begin{align*}
\Gamma & \hookrightarrow G(\QQ_\ell) \cong GL_2(\QQ_\ell), \\
\Gamma & \hookrightarrow G(\ZZ_p) = \MS_2
\end{align*}
induced by the inclusions of the ring $\ZZ[1/\ell]$ into $\QQ_\ell$ and
$\ZZ_p$.  Theorem~\ref{thm:Gamma} says that 
for $\ell$ chosen as above, $\Gamma$ is dense in $\MS_2$ (respectively
$\td{\MS}_2$ for $p=2$).

\subsection{The kernel of the reduced norm}

Define the affine group scheme $G^1$ to be the kernel of reduced norm
$$ G^1 \rightarrow G \xrightarrow{N} \GG_m. $$
The $R$-points of $G^1$ are
given by
$$ G^1(R) = \{x \in \End(C) \otimes R \: : \: N_R(x) = 1 \}. $$
The $\ZZ_p$-points give the closed subgroup
$$ G^1(\ZZ_p) = \MS l_2 $$
of $\MS_2$.  We warn the reader that this group differs from the group
$\MS^1_2$ of \cite{GHMR} and \cite{BehrensK(2)} in that we have not projected
out the Teichm\"uller lift of $\FF_p^\times$ in $\ZZ_p^\times$.  There is
therefore a short exact sequence
$$ 1 \rightarrow \MS l_2 \rightarrow \MS^1_2 \rightarrow \FF_p^\times
\rightarrow 1. $$
We define $\Gamma^1$ to be the group
$$ \Gamma^1 = G^1(\ZZ[1/\ell]). $$
The group $\Gamma^1$ is dense in 
$\MS l_2$ \cite{BehrensLawson}.  

\subsection{Extending by the Galois group}\label{sec:Galois}

Let $Gal \cong \widehat{\ZZ}$ be the Galois group of 
$\br{\FF}_p$ over $\FF_p$.  It is generated by the $p$th power Frobenius
$$ \sigma: \br{\FF}_p \rightarrow \br{\FF}_p. $$
In this section we will introduce compatible actions of $Gal$ on all of our
endomorphism rings.

Recall that given a scheme $X$ over $\br{\FF}_p$, there are three
different Frobenius morphisms, given by the following diagram
$$
\xymatrix{
X \ar[dr]|{\Frob_p^{rel}} \ar@/^/[rrd]^{\Frob_p^{tot}} \ar@/_/[ddr] \\
& X^{(p)} \ar[r]_{\Frob_p} \ar[d] &
X \ar[d] 
\\
& \spec(\br{\FF}_p) \ar[r]_{\sigma^*} &
\spec(\br{\FF}_p)
}
$$
where the scheme $X^{(p)}$ is
the pullback of $X$ over $\sigma^*$.  
The morphism $\Frob_p^{rel}$ is the \emph{relative Frobenius}, 
and $\Frob_p^{tot}$
is the \emph{total Frobenius}.
If $X = Y \otimes_{\FF_p}
\br{\FF}_p$, for a scheme $Y$ over $\FF_p$, then there is a canonical
isomorphism $X \cong X^{(p)}$.  In this case, $\Frob_p$ is an
automorphism of $X$ that covers the automorphism $\sigma$ of $\br{\FF}_p$.

For each $C' \in X^{ss}$, let $\sigma_*C' \in X^{ss}$ be the target of 
the map $\Frob_p$ whose source is $C'$:
$$ \Frob_p: C' \rightarrow \sigma_*C'. $$
Since the curve $C$ was assumed to be defined over $\FF_p$, we have
$\sigma_*C = C$ and $\Frob_p$ takes the form
$$ \Frob_p: C \rightarrow C. $$

For each $\phi \in \End(C)$, the Frobenius $\sigma$ acts on $\phi$ by
$$ \sigma_*\phi = \Frob_p \phi \Frob_p^{-1} \in \End(C). $$
Now if $\phi$ arises from an isogeny defined over $\FF_{p^r}$, then we have 
$(\sigma_*)^r\phi = \phi$.  We conclude that we get an induced
\emph{continuous} action of $Gal$ on $\End(C)$ by ring homomorphisms.
(To be precise, this really should be regarded as an action of $Gal^{op}$,
since iterates of $\Frob_p$ covers the action of $Gal^{op}$ on
$\spec(\br{\FF}_p)$,
but since $Gal$ is abelian, we will ignore this minor point.)

Define $\End_{/\FF_p}(C)$ to be the completed twisted group ring
\begin{align*}
\End_{/\FF_p}(C) & = \End(C)[[Gal]] \\
& = \lim_{r} \End(C)[Gal(\FF_{p^r}/\FF_p)].
\end{align*}
The ring $\End_{/\FF_p}(C)$ consists of endomorphisms of $C$ which do not
cover the identity on $\br{\FF}_p$.

The automorphism $\Frob_p: C \rightarrow C$ does \emph{not} induce a map on
$\br{\FF}_p$-points, because it is not a morphism of schemes over
$\br{\FF}_p$.  The relative Frobenius is not an automorphism of schemes
(since it
is not invertible), 
but it does 
induce an automorphism on $\br{\FF}_p$-points
$$ \Frob_p^{rel}: C(\br{\FF}_p) \rightarrow C(\br{\FF}_p). $$
The following lemma is easily proven using local coordinates.

\begin{lem}\label{lem:FrobFpbar}
Let $\phi$ be an element of $\End(C)$.  Then, on $\br{\FF}_p$-points, 
the endomorphism $\sigma_*\phi$ is given by the composite
$$ \sigma_*\phi : C(\br{\FF}_p) \xrightarrow{(\Frob_p^{rel})^{-1}}
C(\br{\FF}_p) \xrightarrow{\phi} C(\br{\FF}_p) \xrightarrow{\Frob_p^{rel}}
C(\br{\FF}_p). $$
\end{lem}

Let $\sigma^\ZZ$ be the dense subgroup of $Gal$ generated by $\sigma$.
The action of $\sigma^\ZZ$ on $\End(C)$ by ring automorphisms 
induces an action of $\sigma^\ZZ$ on $\Gamma$ by group automorphisms.
Since the norm map is invariant under this action, the action of
$Gal$ restricts to the norm $1$ subgroup $\Gamma^1$.
These actions give rise to extensions
\begin{align*}
\Gamma_{Gal} & = \Gamma \rtimes \sigma^\ZZ, \\
\Gamma^1_{Gal} & = \Gamma^1 \rtimes \sigma^\ZZ.
\end{align*}
We have containments
$$ \Gamma^1_{Gal} \subset \Gamma_{Gal} \subset 
(\End_{/\FF_p}(C)[1/\ell])^\times. $$
In the last containment, the element $\sigma$ of $Gal$ gets mapped to the
automorphism $\Frob_p \in \End_{/\FF_p}(C)$.

The Tate module $T_\ell(C) = \lim_{n} C[\ell^n]$ 
inherits a Galois action through the action of
$\Frob_p^{rel}$, and conjugation by $\Frob_p^{rel}$ induces a Galois action
on $\End_{\ZZ_\ell}(T_\ell(C))$.  Lemma~\ref{lem:FrobFpbar} implies the
following corollary.

\begin{cor}\label{cor:FrobFpbar}
The natural map
$$ \End(C) \rightarrow \End_{\ZZ_\ell}(T_\ell(C)) $$
is Galois equivariant.
\end{cor}

In a manner completely analogous to the case of $\End(C)$, we may define an
action of $Gal$ on $\End(\widehat{C})$ by conjugation with the automorphism 
$\Frob_p$.
The extended Morava stabilizer group is defined by this action:
$$ \GG_2 = \MS_2 \rtimes Gal = \aut_{/\FF_p}(\widehat{C}). $$
The following lemma is clear.

\begin{lem}
The natural map
$$ \End(C) \rightarrow \End(\widehat{C}) $$
is Galois equivariant.
\end{lem}

Thus the inclusion $\Gamma \hookrightarrow \MS_2$ extends to an inclusion 
$$ \Gamma_{Gal} \hookrightarrow \GG_2. $$
Theorem~\ref{thm:Gamma} implies the following proposition.

\begin{prop}
The group $\Gamma_{Gal}$ is dense in $\GG_2$ (respectively, an index $2$
subgroup if $p = 2$).
\end{prop}

\section{The building for $GL_2(\QQ_\ell)$}\label{sec:J'}

\subsection{Construction using lattices}

Let $V$ be a $\QQ_\ell$ vector space of dimension $2$.  The 
Bruhat-Tits building for $GL(V)$ is a contractible $2$-dimensional 
simplicial complex $\mathcal{J}'(V)$ on 
which $GL(V)$ naturally
acts \cite{BruhatTits}.

A lattice $L$ of $V$ is a rank $2$ free $\ZZ_\ell$-submodule such that the
$V = \QQ \otimes L$.   
The complex $\mathcal{J}' = \mathcal{J}'(V)$
is the geometric realization of a semisimplicial
set
$$ \mathcal{J}'_0 \Leftarrow \mathcal{J}'_1 \Lleftarrow \mathcal{J}'_2 $$
where the sets $\mathcal{J}'_i$ are given as the following sets of flags of
lattices in $V$.
\begin{align*}
\mathcal{J}'_0 & = \{ L_0 \: : \: \text{$L_0$ a lattice in $V$}\}, \\
\mathcal{J}'_1 & = \{ L_0 < L_1 \: : \: L_1/L_0
\cong \ZZ/\ell \: \text{or} \: \ZZ/\ell \times \ZZ/\ell \}, \\
\mathcal{J}'_2 & = \{ L_0 < L_1 < L_2 \: : \: 
L_1/L_0 \cong \ZZ/\ell \: \text{and} \: 
L_2/L_0 \cong \ZZ/\ell \times \ZZ/\ell \}. \\
\end{align*}
The $i$th face maps are given by deleting the $i$th terms of the flags.
This semisimplicial set is $GL(V)$ equivariant with the group acting by 
permuting the flags.

We give a description of the underlying topological space of
$\mathcal{J}'$.  Let $\mathcal{J}$ be the $\ell+1$-regular tree (the
building for $SL(V)$).  It is the
infinite tree where every vertex has valence $\ell+1$.

\begin{prop}[See, for instance, \cite{Brown}]
The building $\mathcal{J}'$ is homeomorphic to $\mathcal{J} \times \RR$.
In particular, it is contractible.
\end{prop}

\subsection{Lattices and virtual subgroups}

We now fix $V$ to be the $2$-dimensional $\QQ_\ell$-vector space 
$$ V = V_\ell(C) = T_\ell(C) \otimes \QQ. $$
We shall give a different perspective of the building 
$\mathcal{J}'(V_\ell(C))$
which is more convenient for
understanding the action of the subgroup $\Gamma$ of $GL(V_\ell(C))$.
The lattices of the previous section will be replaced with
certain generalized subgroups of the group $C[\ell^\infty]$, 
which we shall refer to as
virtual subgroups.  These should be thought of as the kernels of
certain quasi-isogenies.  
In this section we define the set of virtual subgroups of $C[\ell^\infty]$,
and show that they are in bijective correspondence with the set of lattices
in $V_\ell(C)$.

Let $\Sub_\ell(C)$ be the set of finite subgroups of $C[\ell^\infty]$.
The set $\Sub_\ell(C)$ carries a natural action of the monoid
$\End_\ell(C)$ (Section~\ref{sec:Gamma}).  Namely, given an endomorphism
$\phi$ in $\End_\ell(C)$, let $\phi$ act on $\Sub_\ell(C)$ by 
$$ \phi: H \mapsto \widehat{\phi}^{-1}(H). $$
The inverse image $\widehat{\phi}^{-1}(H)$ is again a finite subgroup of
$C[\ell^\infty]$. The order of the kernel of $\widehat{\phi}$ 
is the degree
$$\deg(\widehat{\phi}) = \deg (\phi) $$
which is a power of $\ell$.

Observe that the submonoid of $(\ell^k)$th power maps
$$ [\ell^\NN] = \{ [\ell^k] \: : \: k \in \NN \} \subset \End_\ell(C) $$
acts freely on $\Sub_\ell(C)$.
We define the set of \emph{virtual subgroups of $C[\ell^\infty]$} to be the
the set
$$ \Sub^0_\ell(C) = \Sub_\ell(C)[\ell^{-1}] = 
\{ [\ell^k]\cdot H \: : \: k \in \ZZ, H \in \Sub_\ell(C) \} $$
where we have inverted the action of $[\ell]$.  Because the group $\Gamma$ is
given by inverting $\ell$ in the monoid $\End_\ell(C)$
(Lemma~\ref{lem:gpcomplete}), we see that the action of 
$\End_\ell(C)$ on $\Sub^0_{\ell}(C)$
extends to an action of the group $\Gamma$.
Explicitly, for a quasi-isogeny $\psi = \ell^k \phi$ and a virtual subgroup
$H = \ell^{k'} \td{H}$, where $\phi \in \End_\ell(C)$ and 
$\td{H} \in \Sub_\ell(C)$,
this action is given by
$$ \psi \cdot H = \widehat{\psi}^{-1}(H) := 
[\ell^{k+k'}]\cdot \widehat{\phi}^{-1}(\td{H}). $$

As described in Section~\ref{sec:Galois}, the 
set $\Sub_\ell(C)$ possesses a natural Galois action through the
action of $\Frob_p^{rel}$ on $C[\ell^\infty]$ 
$$ \sigma \cdot H = \Frob_p^{rel}(H) $$
and this action extends to $\Sub^0_\ell(C)$.
By Lemma~\ref{lem:FrobFpbar}, this action is compatible with the
Galois action on $\Gamma$, giving $\Sub^0_\ell(C)$ the structure of a
$\Gamma_{Gal}$-set.

The cardinality map
$$ \ord : \Sub_\ell(C) \rightarrow \ell^\NN, $$
which takes a subgroup to its order, 
extends to a map
$$ \ord : \Sub^0_\ell(C) \rightarrow \ell^\ZZ. $$
This map is $\Gamma$-equivariant, where $\phi \in \Gamma$ acts on the 
right-hand
side by multiplication by the degree $N(\phi)$.

Let $\mathcal{L}(V_\ell(C))$ be the set of lattices in $V_\ell(C)$.  It is
a $\Gamma$-set under the inclusion $\Gamma \hookrightarrow GL(V_\ell(C))$,
and the compatible Galois action (Section~\ref{sec:Galois}) induces a
$\Gamma_{Gal}$-action on $\mathcal{L}(V_\ell(C))$.

\begin{prop}\label{prop:kappa}
There is a $\Gamma_{Gal}$-equivariant isomorphism
$$ \kappa: \mathcal{L}(V_\ell(C)) \rightarrow \Sub^0_\ell(C). $$
\end{prop}

\begin{proof}
The map $\kappa$ is the composite
$$ \kappa: \mathcal{L}(V_\ell(C)) \xrightarrow{\td{e}} 
\mathcal{L}(C[\ell^\infty]^* \otimes \QQ) 
\xrightarrow{\ker} \Sub^0_\ell(C) $$
where $\td{e}$ is the Galois equivariant isomorphism of 
Proposition~\ref{prop:Weil}, under
which $C[\ell^\infty]^*$ inherits the
$\Gamma$ action given in Lemma~\ref{lem:Weilaction}.  
We describe the map $\ker$.  If we are given a lattice $L \subset
C[\ell^\infty]^* \otimes \QQ$ which is
actually contained in 
$C[\ell^\infty]^*$, 
we define $\ker(L)$ to be
the subgroup of $C$ given by
$$ \ker(L) = \bigcap_{\alpha \in L} \ker{\alpha} \in \Sub_\ell(C). $$
If $L$ is a general lattice in 
$C[\ell^\infty]^*\otimes \QQ$, then there exists a $k$ such that $\ell^k L$ is
contained in $C[\ell^\infty]^*$.  Define $\ker(L)$ to be the virtual subgroup 
$$ \ker(L) = [\ell^{-k}] \cdot \ker(\ell^k L). $$
This is easily seen to be independent of the choice of $k$.

In order to show that $\ker$ is $\Gamma$-equivariant,
it suffices to check the $\Gamma$-equivariance on elements $\phi$ of $\Gamma$
contained in $\End_\ell(C)$ and on lattices contained in
$C[\ell^\infty]^*$.
We then have
\begin{align*}
\ker(\phi \cdot L) & = \ker(\{ \alpha \circ \widehat{\phi} \: : \: 
\alpha \in L\}) \\
& = \bigcap_{\alpha \in L} \ker(\alpha \circ \widehat{\phi}) \\
& = \bigcap_{\alpha \in L} \widehat{\phi}^{-1}( \ker(\alpha)) \\
& = \widehat{\phi}^{-1}(\ker(L)) \\
& = \phi \cdot \ker(L).
\end{align*}
The Galois equivariance of $\ker$ is easily verified:
\begin{align*}
\ker(\sigma \cdot L)
& = \bigcap_{\alpha \in L} \ker(\Frob_p^{rel} \circ \alpha \circ
(\Frob_p^{rel})^{-1}) \\
& = \bigcap_{\alpha \in L} \ker(\alpha \circ
(\Frob_p^{rel})^{-1}) \\
& = \bigcap_{\alpha \in L} \Frob_p^{rel}(\ker \alpha) \\
& = \Frob_p^{rel} (\ker(L)) \\
& = \sigma \cdot \ker(L).
\end{align*}

To see that the map $\ker$ is a bijection, we construct an inverse.  Given
a finite subgroup $H$ of $C[\ell^\infty]$, define $L_H$ to be the subgroup 
of $C[\ell^\infty]^*$ given by
$$ L_H = \{ \alpha \in C[\ell^\infty]^* \: : \: \alpha(H) = 1 \}. $$
We claim that $L_H$ is a lattice, that is, that the inclusion
$$ L_H \otimes \QQ \hookrightarrow C[\ell^\infty]^*\otimes \QQ $$
is a bijection.  It suffices to show that there exists a positive
integer $k$ such that $C[\ell^\infty]^*$ is contained in $\ell^{-k} L_H$, or
equivalently, such that $\ell^k C[\ell^\infty]^*$ is contained in 
$L_H$.  This is
accomplished by choosing $k$ such that $H$ is contained in the $\ell^k$
torsion of $C$.  The mapping 
$ H \mapsto L_H $
extends to a map
$$ \Sub^0_\ell(C) \rightarrow \mathcal{L}(C[\ell^\infty]^* \otimes
\QQ) $$
which is inverse to the map $\ker$.
\end{proof}

Suppose that $H = [\ell^k] \cdot \td{H}$ and $H' = [\ell^{k'}]\cdot
\td{H}'$ 
are
virtual subgroups of $C[\ell^\infty]$, for $\td{H},\td{H}' \in \Sub_\ell(C)$.  
We may assume that $k = k'$.  We
shall say that $H$ is \emph{contained} in $H'$ and write 
$$ H \le H' $$
if $\td{H}$ is contained in $\td{H}'$.  Define the \emph{quotient} to be
$$ H'/H = \td{H}'/\td{H}. $$
Observe that this depends on the choice of $k$, but any two choices of $k$
will give canonically isomorphic quotients.

\begin{lem}\label{lem:quotients}
Suppose that $L_0 \le L_1$ are lattices in $V_\ell(C)$.  Then there is a
containment
$$ \kappa(L_1) \le \kappa(L_0) $$
of virtual subgroups and a (non-canonical) isomorphism
$$ L_1/L_0 \cong \kappa(L_0)/\kappa(L_1)$$
between the quotients.
\end{lem}

\begin{proof}
Let the lattices $\td{e}(L_i) \subset C[\ell^\infty]^* \otimes \QQ$ be the 
images
of the lattices $L_i$ under the map $\td{e}$ of
Proposition~\ref{prop:Weil}.
We may as well assume that the lattices $\td{e}(L_i)$ are contained in 
$C[\ell^\infty]^*$.  The subgroups $\kappa(L_i)$ are
the kernels of the dual projections
$$
0 \rightarrow \kappa(L_i) \rightarrow C[\ell^\infty] \rightarrow
\td{e}(L_i)^*
\rightarrow 0.
$$
There are therefore isomorphisms
$$ (L_1/L_0)^* \cong \ker(L_1^* \twoheadrightarrow L_0^*) \cong
\kappa(L_0)/\kappa(L_1) $$
using the exactness of the Pontryagin dual and the $3 \times 3$ lemma.
Since $L_1/L_0$ is finite, it is non-canonically isomorphic to its
Pontryagin dual $(L_1/L_0)^*$.
\end{proof}

\subsection{Construction of $\mathcal{J}'$ using virtual
subgroups}\label{sec:J'virtual}

The map $\kappa$ of Proposition~\ref{prop:kappa} and
Lemma~\ref{lem:quotients}
gives the following alternative description of 
the sets of simplices in the semi-simplicial set
$$ \mathcal{J}'_0 \Leftarrow \mathcal{J}'_1 \Lleftarrow \mathcal{J}'_2 $$
in terms of flags of virtual subgroups of $C[\ell^\infty]$:
\begin{align*}
\mathcal{J}'_0 & = \{ H_0 \: : \: \text{$H_0$ a virtual subgroup 
of $C[\ell^\infty]$}\}, \\
\mathcal{J}'_1 & = \{ H_1 < H_0 \: : \: H_0/H_1
\cong \ZZ/\ell \: \text{or} \: \ZZ/\ell \times \ZZ/\ell \}, \\
\mathcal{J}'_2 & = \{ H_2 < H_1 < H_0 \: : \: 
H_1/H_2 \cong \ZZ/\ell \: \text{and} \: 
H_0/H_2 \cong \ZZ/\ell \times \ZZ/\ell \}. \\
\end{align*}
The $i$th face maps are given by deleting the $i$th terms of the flags.
This semisimplicial set is $\Gamma$ equivariant with the group acting by 
permuting the flags.  This action agrees with the action given by the
inclusion $\Gamma \hookrightarrow GL(V_\ell(C))$ since the map $\kappa$ was
proven to be $\Gamma$-equivariant.

\subsection{The $\Gamma$ orbits in $\mathcal{J}'$}\label{sec:orbits}

We shall explicitly 
identify the $\Gamma$ orbits of the sets $\mathcal{J}'_i$, and determine
their isotropy.

Recall that $X^{ss}$ is the set of isomorphism classes of supersingular
elliptic curves $C'$
defined over $\br{\FF}_p$, and $X^{ss}_0(\ell)$ is the set of
isomorphism classes of pairs $(C',H)$ of supersingular curves $C'$ with a
cyclic subgroup $H$ of order $\ell$.  Fix representatives of these
isomorphism classes.  We shall make use of the following result of Kohel.

\begin{thm}[Kohel, {\cite[Cor.~77]{Kohel}}]\label{thm:Kohel}
Let $C'$ and $C''$ be supersingular elliptic curves over $\br{\FF}_{p}$.
Then for all $k \gg 0$, there exists an isogeny $\phi: C' \rightarrow C''$
of degree $\ell^k$.
\end{thm}

Since there are finitely many elliptic curves, there exists an $e > 0$ so
that we may choose isogenies
$$ \phi_{C'} : C \rightarrow C' $$
of degree $\ell^{2e}$ for every $C' \in X^{ss}$.
We may as well assume that $\phi_C = [\ell^e]$.  
This uniformity in the degrees of the isogenies $\phi_{C'}$ has the effect
of simplifying some of our proofs.  
Our insistence on using isogenies of degree an even power of
$\ell$ will come into play in Section~\ref{sec:orbits1}
(see Lemma~\ref{lem:isogenyparity}).

Our choices of the isogenies $\phi_{C'}$ give embeddings 
$$ \iota_{C'} : \aut(C') \hookrightarrow \Gamma $$
of the
automorphism groups $\aut(C')$ for $C' \in X^{ss}$.  Given an automorphism
$\gamma \in \aut(C')$, let $\iota_{C'}(\gamma)$ be the quasi-isogeny of $C$
given by 
$$ \iota_{C'}(\gamma) = [\ell^{-2e}]\cdot (\widehat{\phi}_{C'} \circ 
\gamma \circ
\phi_{C'}). $$
The factor of $[\ell^{-2e}]$ makes $\iota_{C'}$ a homomorphism of
groups.  
For $(C',H) \in X^{ss}_0(\ell)$ we regard the subgroups 
$\aut(C',H)$ of $\aut(C')$ to be
embedded in $\Gamma$ by $\iota_{C'}$.

\begin{prop}\label{prop:orbits}
The $\Gamma$-sets $\mathcal{J}'_i$ decompose into $\Gamma$-orbits as
follows:
\begin{align*}
\mathcal{J}'_0 & = \coprod_{C' \in X^{ss}} \mathcal{J}'_0[C'], \\
\mathcal{J}'_1 & = \coprod_{(C',H) \in X^{ss}_0(\ell) } \mathcal{J}'_1[C',H]
\quad \amalg \quad \coprod_{C' \in X^{ss}} \mathcal{J}'_1[C'], \\
\mathcal{J}'_2 & = \coprod_{(C',H) \in X^{ss}_0(\ell)}
\mathcal{J}'_2[C',H].
\end{align*}
These orbits are given as follows:
\begin{align*}
\mathcal{J}'_0[C'] & = \{ H_0 \: : \: C/\td{H}_0 \cong C' \}, \\
\mathcal{J}'_1[C',H] & = \{ H_1 < H_0 \: : \: (C/\td{H}_0, 
\ell \cdot \td{H}_1/\td{H}_0) 
\cong (C',H) \}, \\
\mathcal{J}'_1[C'] & = \{ H_1 < H_0 \: : \:  C/\td{H}_0 \cong
C' \: \text{and} \: H_0 = \ell \cdot H_1 \}, \\
\mathcal{J}'_2[C',H] & = \{ H_2 < H_1 < H_0 \: : \:  (C/\td{H}_0,
\ell \cdot \td{H}_1/\td{H}_0) \cong (C',H) \\
& \quad \quad \quad \quad \quad \quad \quad \quad \quad \quad \quad
\text{and} \: H_0 = \ell \cdot H_2 \}
\end{align*}
where the subgroups $\td{H}_i$ of $C$ are obtained from the virtual subgroups
$H_i$ by multiplying by a suitable uniform power of $\ell$. 
\end{prop}

\begin{proof}
We explain the decomposition of $\mathcal{J}'_0$.  The arguments are
essentially the same for the other $\mathcal{J}'_i$.  
We first verify that the set $\mathcal{J}'_0[C']$ is closed under the 
action of $\Gamma$.  It suffices to check that if $H_0$ is a subgroup of
$C$ such that $C/H_0 \cong C'$ and $\phi$ is an isogeny in $\End_\ell(C)$,
then $\phi \cdot H_0 = \widehat{\phi}^{-1}(H_0)$ has the property that
there is an isomorphism $C/\widehat{\phi}^{-1}(H_0) \cong C'$.  
This isomorphism
exists because $\widehat{\phi}^{-1}(H_0)$ is the kernel of
the composite
$$ C \xrightarrow{\widehat{\phi}} C \rightarrow C/H_0 \xrightarrow{\cong}
C'. $$
We now verify transitivity.  Suppose that $H_0$ and $H_0'$ are 
elements of
$\mathcal{J}'_0[C']$.  We may assume that the virtual subgroups
$H_i$ are actually subgroups.  Let $\phi_{{H}_0}$ and
$\phi_{{H}_0'}$ be the 
quotient maps
\begin{gather*}
\phi_{{H}_0}: C \rightarrow C/{H}_0, \\
\phi_{{H}_0'}: C \rightarrow C/{H}_0'. \\
\end{gather*}
By hypothesis, there is an isomorphism
$$ \gamma: C/{H}_0 \xrightarrow{\cong} C/{H}_0'.$$
Let $\psi$ be the composite
$$ \psi: C \xrightarrow{\phi_{{H}_0}} C/{H}_0 \xrightarrow{\gamma} C/{H}_0' 
\xrightarrow{\widehat{\phi}_{{H}_0'}} C. $$ 
Then we have
$$ \psi \cdot {H}_0 = [\ell^i] \cdot {H}_0' $$
where $\ell^i$ is the order of $H_0$, 
so $(\ell^{-i}\psi) \cdot {H}_0 = H_0'$.  The action of $\Gamma$ is
therefore transitive.
\end{proof}

\begin{prop}\label{prop:isotropy}
The orbits of Proposition~\ref{prop:orbits} are identified as follows:
\begin{align*}
\mathcal{J}'_i[C'] & \cong \Gamma/\aut(C'), \\
\mathcal{J}'_i[C',H] & \cong \Gamma/\aut(C',H).
\end{align*}
\end{prop}

\begin{proof}
For each $C'$ in $X^{ss}$, let $K_{C'}$ be the kernel of the isogeny
$\phi_{C'}$.  For each $(C',H)$ in $X^{ss}_0(\ell)$, let $K_{(C',H)}$
be the kernel of the composite
$$ C \xrightarrow{\phi_{C'}} C' \twoheadrightarrow C'/H. $$
Then the stabilizers of certain well-chosen points of $\mathcal{J}'_i$ are
easily determined:
\begin{align*}
& \Sta_{\Gamma}(K_{C'}) = \aut(C'), \\
& \Sta_{\Gamma}(\ell^{-1} \cdot K_{(C',H')} < K_{C'}) = \aut(C',H), \\
& \Sta_{\Gamma}(\ell^{-1} \cdot K_{C'} < K_{C'}) = \aut(C'), \\
& \Sta_{\Gamma}(\ell^{-1} \cdot K_{C'} < \ell^{-1} \cdot K_{(C',H)} 
< K_{C'}) = \aut(C',H).
\end{align*}
\end{proof}

Combining Propositions~\ref{prop:orbits} and \ref{prop:isotropy}, we have
the following theorem.

\begin{thm}\label{thm:calJ'}
There are $\Gamma_{Gal}$-equivariant isomorphisms
\begin{align}
\begin{split}
\mathcal{J}'_0 & \cong \coprod_{C' \in X^{ss}} \Gamma/\aut(C'), \\
\mathcal{J}'_1 & \cong \coprod_{(C',H) \in X^{ss}_0(\ell) }
\Gamma/\aut(C',H)
\quad \amalg \quad \coprod_{C' \in X^{ss}} \Gamma/\aut(C'), \\
\mathcal{J}'_2 & \cong \coprod_{(C',H) \in X^{ss}_0(\ell)}
\Gamma/\aut(C',H).
\end{split}\label{eq:J'gamma}
\end{align}
\end{thm}

\begin{proof}
We are only left with verifying the Galois equivariance of these
isomorphisms.  We will only treat the case of $\mathcal{J}'_i$ with $i = 0$.  
The cases of $i > 0$ are completely analogous.

Recall that we have embedded $\aut(C')$ in $\Gamma$ by conjugating with
$\phi_{C'}$.  Our practice of denoting the $\Gamma$ orbit corresponding to $C'$
by $\Gamma/\aut(C')$ is misleading, because it conceals the manner in which
we have embedded $\aut(C')$.  The orbit is more precisely
given by
$$ \Gamma/\phi_{C'}^{-1} \aut(C') \phi_{C'}. $$
The decomposition of $\calJ'_0$ is given by the composite
$$ \Gamma/\phi_{C'}^{-1} \aut(C') \phi_{C'} \xrightarrow[f]{\cong} 
\Gamma \cdot K_{C'} \xrightarrow{\cong} \calJ'_0[C']
$$
where $K_{C'} \in \Sub^0_\ell(C)$ is the kernel of $\phi_{C'}$.  The map
$f$ is given by
$$ f(x\phi_{C'}^{-1} \aut(C') \phi_{C'}) = x \cdot K_{C'} $$
for $x \in \Gamma$.

The action of $\sigma \in Gal$
$$ \sigma : \calJ'[C'] \rightarrow \calJ'[\sigma_*C'] $$
is given by
$$ \sigma \cdot H = \Frob_p^{rel}(H) $$
for $H \in \Sub_\ell(C)$ with $C/H \cong C'$.
Using Lemma~\ref{lem:FrobFpbar}, we have
$$ \Frob_p^{rel}(K_{C'}) = y_{C'} \cdot K_{\sigma_*C'} $$
for
$$ y_{C'} = \ell^{-2e} \cdot (\sigma_* \widehat{\phi}_{C'}) \cdot
\phi_{\sigma_*C'} \in \Gamma. $$
Thus the compatible $\sigma$ action
$$ \sigma: \Gamma \cdot K_{C'} \rightarrow \Gamma \cdot K_{\sigma_*C'} $$
is given by
\begin{align*}
\sigma \cdot (x \cdot K_{C'})
& = (\sigma_*x) \cdot \sigma \cdot K_{C'} \\
& = (\sigma_*x) \cdot y_{C'} \cdot K_{\sigma_*C'}.
\end{align*}

We now compute the image of the subgroup $\phi_{C'}^{-1}\aut(C')\phi_{C'}$
under the action of $\sigma_*$ on $\Gamma$.
\begin{align*}
\sigma_*(\phi_{C'}^{-1} \aut(C') \phi_{C'}) 
& = \Frob_p\phi_{C'}^{-1} \aut(C') \phi_{C'} \Frob_p^{-1} \\
& = (\sigma_*\phi_{C'})^{-1} \Frob_p\aut(C')\Frob_p^{-1}
(\sigma_*\phi_{C'}) \\
& = (\sigma_*\phi_{C'})^{-1} \aut(\sigma_*C')
(\sigma_*\phi_{C'}) \\
& = y_{C'} \phi_{\sigma_*C'}^{-1} \aut(\sigma_*C') \phi_{\sigma_*C'}
y_{C'}^{-1}. 
\end{align*}
With this in mind, the natural action of $\sigma_*$ on 
$\Gamma/\phi_{C'}^{-1}\aut(C')\phi_{C'}$ is
given by
\begin{align}
\begin{split}\label{eq:GalJ'}
\sigma: \Gamma/\phi_{C'}^{-1}\aut(C')\phi_{C'} & \rightarrow 
\Gamma/\phi_{\sigma_*C'}^{-1}\aut(\sigma_*C')\phi_{\sigma_*C'} \\
x \cdot \phi_{C'}^{-1}\aut(C')\phi_{C'} & \mapsto 
\sigma_*x \cdot y_{C'} \cdot 
\phi_{\sigma_*C'}^{-1}\aut(\sigma_*C')\phi_{\sigma_*C'}.
\end{split}
\end{align}
This action makes the map $f$ Galois equivariant.
\end{proof}

\subsection{The semi-simplicial structure of $\calJ'$}\label{sec:J'face}

In Section~\ref{sec:orbits} we gave a $\Gamma$-equivariant orbit
decomposition of the $n$-simplices of $\calJ'_\bullet$.  In this section we
shall describe the face maps in terms of this orbit decomposition.

We shall first need some definitions.
Recall from Section~\ref{sec:orbits} that we have fixed representatives
$C'$ (respectively $(C',H)$) for each isomorphism class of $X^{ss}$
(respectively $X^{ss}_0(\ell)$).  We also fixed isogenies
$$ \phi_{C'} : C \rightarrow C' $$
of degree $\ell^{2e}$ for each $C' \in X^{ss}$.

For each pair $(C',H) \in X^{ss}_0(\ell)$, there is an induced degree $\ell$ 
isogeny given by the quotient
$$ q_H : C' \rightarrow C'/H. $$
Let $\widehat{H} \subset C'/H$ be the kernel of the dual isogeny
$\widehat{q}_H$.  Then there exists a pair $(C_H, d(H)) \in
X_0^{ss}(\ell)$ and an isomorphism
$$ \alpha_H: C'/H \rightarrow C_H $$
which sends $\widehat{H}$ to $d(H)$.  Define $\phi_{H}$ to be the composite
$$ \phi_H: C' \xrightarrow{q_H} C'/H \xrightarrow{\alpha_H} C_H. $$
Then $d(H)$ is the kernel of the dual isogeny $\widehat{\phi}_{H}$.

For each pair $(C',H) \in X^{ss}_0(\ell)$, define elements of $g_{(C',H)}$
of $\Gamma$ by
\begin{equation}\label{eq:g_C'H}
g_{(C',H)} = \ell^{-2e-1} \cdot (\widehat{\phi}_{C'} \circ 
\widehat{\phi}_H \circ \phi_{C_H}). 
\end{equation}

\begin{prop}\label{prop:J'face}
Under the isomorphisms of Equation~(\ref{eq:J'gamma}), the face maps of the
semisimplicial set $\calJ'_\bullet$ are given as follows.
\begin{align*}
d_i: & \calJ'_1 \rightarrow \calJ'_0 \\
& d_0(x\aut(C',H)) = xg_{(C',H)}\aut(C_H), \\
& d_1(x\aut(C',H)) = x\aut(C'), \\
& d_0(x\aut(C')) = x\cdot \ell^{-1}\aut(C'), \\
& d_1(x\aut(C')) = x\aut(C'). \\
d_i: & \calJ'_2 \rightarrow \calJ'_1 \\
& d_0(x\aut(C',H)) = xg_{(C',H)}\aut(C_H,d(H)), \\
& d_1(x\aut(C',H)) = x\aut(C'), \\
& d_2(x\aut(C',H)) = x\aut(C',H). \\
\end{align*}
\end{prop}

\begin{proof}
We simply must evaluate the face maps of $\calJ'$, as given in
Section~\ref{sec:J'virtual} on the orbit representatives chosen
in the proof of Proposition~\ref{prop:isotropy}:
\begin{align*}
d_0(\ell^{-1}K_{(C',H)} < K_{C'}) & = \ell^{-1}K_{(C',H)} \\
& = g_{(C',H)} \cdot K_{C_H}, \\
d_1(\ell^{-1}K_{(C',H)} < K_{C'}) & = K_{C'}, \\
d_0(\ell^{-1}K_{C'} < K_{C'}) & = \ell^{-1}K_{C'}, \\
d_1(\ell^{-1}K_{C'} < K_{C'}) & = K_{C'}, \\
d_0(\ell^{-1}K_{C'} < \ell^{-1}K_{(C',H)} < K_{C'}) & = 
\ell^{-1}K_{C'} < \ell^{-1}K_{(C',H)} \\
& = g_{(C',H)} \cdot (\ell^{-1}K_{(C_H,d(H))} < K_{C_H}), \\
d_1(\ell^{-1}K_{C'} < \ell^{-1}K_{(C',H)} < K_{C'}) & =
\ell^{-1}K_{C'} < K_{C'}, \\
d_2(\ell^{-1}K_{C'} < \ell^{-1}K_{(C',H)} < K_{C'}) & =
\ell^{-1}K_{(C',H)} < K_{C'}.
\end{align*}
\end{proof}

\section{The building for $SL_2(\QQ_\ell)$}\label{sec:J}

In this section we recall the construction of 
the building $\mathcal{J}$ for $SL_2(\QQ_\ell)$.  We then give a
reinterpretation in terms of virtual subgroups of $C[\ell^\infty]$ which is
more amenable to understanding the action of the subgroup $\Gamma^1$
of $SL(V_\ell(C))$.

We decompose $\mathcal{J}$ into orbits under the action of $\Gamma^1$, and 
demonstrate that this group acts on $\calJ$ with finite
stabilizers.  We then explain how
Bass-Serre theory gives the structure of the group $\Gamma^1$
as the fundamental groups of a graph of 
finite groups.    
Much of the material in this section may be found in
Serre's book \cite{Serre}.

\subsection{The construction of $\mathcal{J}$ using lattices}

Let $V$ be a $\QQ_\ell$ vector space of dimension $2$.  Two lattices $L$ and
$L'$ in $V$ are said to be \emph{homothetic} if there exists a $c \in
\QQ_\ell^\times$ such that 
$$ L' = cL. $$
Since $L$ and $L'$ are $\ZZ_\ell$-modules, $c$ may be chosen to be $\ell^k$
for some integer $k$.
The construction of $\calJ$ follows the construction of $\calJ'$ except
that we use homothety classes of lattices in $V$.  The building $\calJ$
is a $1$-dimensional contractible simplicial complex on which $SL(V)$ acts.
Topologically, $\calJ$ is an $\ell+1$-regular tree.

Specifically, $\calJ$ is the geometric realization of a semi-simplicial
$SL(V)$ set of the form
$$ \calJ_0 \leftleftarrows \calJ_1 $$
where the sets $\calJ_i$ are sets of flags of homothety classes of
lattices in $V$:
\begin{align*}
\mathcal{J}_0 & = \{ [L_0] \: : \: \text{$[L_0]$ a homothety class of 
lattice in $V$}\}, \\
\mathcal{J}_1 & = \{ \{ [L_0],[L_1] \} \: : \: 
\text{there exist reps $L_0 < L_1$ such that
$L_1/L_0
\cong \ZZ/\ell$} \}.
\end{align*}
The group $GL(V)$ acts by permuting the lattice classes in the flags.  This
action restricts to an action of $SL(V)$. Since we are taking homothety
classes of lattices, the center $\QQ_\ell^\times
\subseteq GL(V)$ acts trivially on $\calJ$, so the action also factors through
$PGL(V)$.

There is a $GL(V)$-equivariant projection
$$ \nu: \calJ' \rightarrow \calJ $$
given by taking homothety classes of the lattices that make up the
flags of $\calJ'$.  Under this map, the simplices of $\calJ'_1$
corresponding to flags $L_0 < L_1$ with $L_1/L_0 \cong \ZZ/\ell \times
\ZZ/\ell$, as well as all of the simplices of $\calJ'_2$, become
degenerate. 

\subsection{The construction of $\mathcal{J}$ using virtual
subgroups}\label{sec:Jvirtual}

Let $V = V_\ell(C)$.
The same methods that construct $\calJ'$ in terms of virtual subgroups
construct $\calJ$ in terms of homothety classes of virtual subgroups of
$C[\ell^\infty]$.  Here, two virtual subgroups $H$ and $H'$ are
said to be homothetic if there exists an integer $k$ such that 
$$ H' = [\ell^k]\cdot H. $$

\begin{lem}\label{lem:subgroupparity}
Every virtual subgroup  $H$ of $C[\ell^\infty]$ is homothetic to a unique
virtual subgroup  $H'$ where the order of $H'$ is either $1$ or
$\ell$.  
The virtual subgroup  $H'$ is uniquely expressible in the form
$$ H' = [\ell^k] \cdot H'' $$
for some integer $k$, where $H''$ is a subgroup of $C$ isomorphic to
$\ZZ/\ell^m$.
\end{lem}

\begin{proof}
The virtual subgroup $H'$ is $[\ell^{-i}]\cdot H$ where the order of $H$
is either $\ell^{2i}$ or $\ell^{2i+1}$.  To produce the canonical
representative $H''$, we may as well assume that the representative $H$ of the
homothety class $[H]$ is a subgroup.  Let $j$ be maximal so that the
$\ell^j$-torsion subgroup $C[\ell^j]$ is contained in $H$.  Then the
subgroup $H''$ is given by
$$ H'' = H/C[\ell^j] \subset C/C[\ell^j]. $$
We may regard $H''$ as being contained in $C$ under the canonical
isomorphism $[\ell^j]: C/C[\ell^j] \xrightarrow{\cong} C$.
\end{proof}

If the order of the group $H'$ given by Lemma~\ref{lem:subgroupparity} is $1$, 
we shall say 
the homothety class $[H]$ is \emph{even}.  Otherwise we shall say that the
homothety class $[H]$ is
\emph{odd}.

Observe that the isomorphism $\kappa$ of Proposition~\ref{prop:kappa}
identifies homothety classes of lattices with homothety classes of
virtual subgroups.
The semi-simplicial set
$$ \calJ_0 \leftleftarrows \calJ_1 $$
whose realization is $\calJ$ may therefore be described in terms of
virtual subgroups:
\begin{align*}
\mathcal{J}_0 & = \{ [H_0] \: : \: \text{$[H_0]$ a homothety class of 
virtual subgroup  in $C[\ell^\infty]$}\}, \\
\mathcal{J}_1 & = \{ ([H_1],[H_0]) \: : \: 
\text{there exist reps $H_1 < H_0$ such that
$H_0/H_1
\cong \ZZ/\ell$,} \\
& \quad\quad\quad\quad\quad\quad\quad\quad\quad
\text{$[H_0]$ even, $[H_1]$ odd}
\}.
\end{align*}
The group $\Gamma^1$  
acts by permuting the classes of virtual subgroups.

\subsection{The $\Gamma^1$ orbit decomposition of 
$\mathcal{J}$}\label{sec:orbits1}

We must first remark that $\Gamma^1$ contains only half of the isogenies of
$\Gamma$ modulo $[\ell^\ZZ]$.

\begin{lem}\label{lem:isogenyparity}
Every quasi-isogeny $\phi \in \Gamma^1$ is expressible uniquely in the form
$$ \phi = \ell^{-i} \phi' $$
where $\phi'$ is an endomorphism of $C$ whose kernel is isomorphic to
$\ZZ/\ell^{2i}$.
\end{lem}

\begin{proof}
There is a short exact sequence
$$ 1 \rightarrow \Gamma^1 \rightarrow \Gamma \xrightarrow{N} \ell^\ZZ
\rightarrow 1. $$
The lemma is immediate from the fact that $\Gamma = \End_\ell(C)[\ell^{-1}]$
and $N([\ell]) = \ell^2$.
\end{proof}

The orbits of $\calJ_i$ are given in the following proposition, whose proof
is completely analogous to that of Proposition~\ref{prop:orbits}.  The
decomposition of $\calJ_0$ into the two parity classes of $\Gamma^1$ orbits 
is a consequence of
Lemmas~\ref{lem:subgroupparity} and \ref{lem:isogenyparity}.

\begin{prop}\label{prop:orbits1}
The $\Gamma^1$-sets $\mathcal{J}_i$ decompose into $\Gamma^1$-orbits as
follows:
\begin{align*}
\mathcal{J}_0 & = \coprod_{C' \in X^{ss}} (\mathcal{J}_0[C']_{\it even}
\amalg \mathcal{J}_0[C']_{\it odd}), \\
\mathcal{J}_1 & = \coprod_{(C',H) \in X^{ss}_0(\ell) }
\mathcal{J}_1[C',H].
\end{align*}
These orbits are given as follows:
\begin{align*}
\mathcal{J}_0[C']_{\it even} & = \{ [H_0] \: : \: C/H_0 \cong C' \:
\text{and $[H_0]$ even} \}, \\
\mathcal{J}_0[C']_{\it odd} & = \{ [H_0] \: : \: C/H_0 \cong C' \:
\text{and $[H_0]$ odd} \}, \\
\mathcal{J}_1[C',H] & = \{ ([H_1],[H_0]) \: : \: 
\text{there exist reps $H_1 < H_0$ such that $H_0/H_1 \cong \ZZ/\ell$,} \\ 
& \quad\quad\quad\quad\quad\quad\quad\quad\quad
\text{$[H_0]$ even, $[H_1]$ odd, $(C/H_0, \ell \cdot H_1/H_0) \cong
(C',H)$} \}.
\end{align*}
\end{prop}

Recall from Section~\ref{sec:orbits} that we have embedded the group
$\aut(C')$ as a subgroup of $\Gamma^1$ by conjugating by the isogeny
$\phi_{C'}$:
\begin{align*}
\iota_{C'}: & \aut(C') \hookrightarrow \Gamma^1 \\
& \alpha \mapsto \phi_{C'}^{-1} \alpha \phi_{C'}.
\end{align*}
Fix an endomorphism $\phi$ of $C$ of degree $\ell^{2r+1}$ for $r \gg 0$.  Such
an endomorphism exists by Theorem~\ref{thm:Kohel}.  We shall use
$\br{\aut(C')}$ to denote the image of the 
different embedding of $\aut(C')$ in
$\Gamma^1$ given by
\begin{align*}
\br{\iota}_{C'}: & \br{\aut(C')} \hookrightarrow \Gamma^1 \\
& \alpha \mapsto \phi^{-1} \phi_{C'}^{-1} \alpha \phi_{C'} \phi.
\end{align*}
The isotropy of $\calJ$ is described in the following proposition.

\begin{prop}\label{prop:isotropy1}
The orbits of Proposition~\ref{prop:orbits1} are given by:
\begin{align*}
\mathcal{J}_0[C']_{\it even} & \cong \Gamma^1/\aut(C'), \\
\mathcal{J}_0[C']_{\it odd} & \cong \Gamma^1/\br{\aut(C')}, \\
\mathcal{J}_1[C',H] & \cong \Gamma^1/\aut(C',H).
\end{align*}
\end{prop}

\begin{proof}
For each $C'$ in $X^{ss}$, let $K_{C'}^{\it even} = K_{C'}$ be the 
kernel of the 
isogeny $\phi_{C'}$, and let $K^{odd}_{C'}$ be the kernel of the composite
$$ C \xrightarrow{\phi} C \xrightarrow{\phi_{C'}} C'. $$
For each $(C',H)$ in $X^{ss}_0(\ell)$, let $K_{(C',H)}$
be the kernel of the composite
$$ C \xrightarrow{\phi_{C'}} C' \twoheadrightarrow C'/H. $$
Then the stabilizers of certain well-chosen points of $\mathcal{J}_0[C']_{\it
even}$, $\calJ_0[C']_{\it odd}$, and $\calJ_1[C',H]$ are
easily determined:
\begin{align*}
& \Sta_{\Gamma^1}([K_{C'}^{\it even}]) = \aut(C'), \\
& \Sta_{\Gamma^1}([K_{C'}^{\it odd}]) = \br{\aut(C')}, \\
& \Sta_{\Gamma^1}([K_{(C',H)}],[K_{C'}^{\it even}]) = \aut(C',H).
\end{align*}
\end{proof}

\subsection{The semi-simplicial structure of $\calJ$}\label{sec:Jface}

In this section we
shall describe the face maps in the semi-simplicial set $\calJ_\bullet$ in
terms of the orbit decomposition given in Section~\ref{sec:orbits1}.

Combining Propositions~\ref{prop:orbits1} and \ref{prop:isotropy1}, we have
$\Gamma^1$-equivariant isomorphisms
\begin{align}
\begin{split}
\mathcal{J}_0 & \cong \coprod_{C' \in X^{ss}} \Gamma^1/\aut(C') \amalg
\Gamma^1/\br{\aut(C')}, \\
\mathcal{J}_1 & \cong \coprod_{(C',H) \in X^{ss}_0(\ell) }
\Gamma^1/\aut(C',H).
\end{split}\label{eq:Jgamma}
\end{align}
For each pair $(C',H) \in X^{ss}_0(\ell)$, define an element $g^1_{(C',H)}$
of $\Gamma^1$ by
$$ g^1_{(C',H)} = \ell^{-(r+2e+1)} \cdot (\widehat{\phi}_{C'} \circ 
\widehat{\phi}_H \circ \phi_{C_H} \circ \phi). $$

\begin{prop}\label{prop:Jface}
Under the isomorphisms of Equation~(\ref{eq:Jgamma}), the face maps of the
semisimplicial set $\calJ_\bullet$ are given as follows:
\begin{align*}
d_i: & \calJ_1 \rightarrow \calJ_0 \\
& d_0(x\aut(C',H)) = xg^1_{(C',H)}\br{\aut(C_H)}, \\
& d_1(x\aut(C',H)) = x\aut(C'). \\
\end{align*}
\end{prop}

\begin{proof}
We evaluate the face maps of $\calJ$, as given in
Section~\ref{sec:Jvirtual} on the orbit representatives chosen
in the proof of Proposition~\ref{prop:isotropy1}:
\begin{align*}
d_0([K_{(C',H)}],[K_{C'}^{\mathit{even}}]) & = [K_{(C',H)}] \\
& = g^1_{(C',H)} \cdot [K_{C_H}^{\mathit{odd}}], \\
d_1([K_{(C',H)}],[K_{C'}^{\mathit{even}}]) & = [K_{C'}^\mathit{even}].
\end{align*}
\end{proof}

\subsection{The structure of $\Gamma^1$}

A graph of groups is a graph $Y$ 
whose vertices and edges are labeled with 
finite groups, with inclusions compatible with the gluing data of the graph.
In \cite{Serre}, the notion of the fundamental group of a connected graph of 
groups is given.  If the graph $Y$ is a tree, then this fundamental group
is simply a suitable amalgamation of the labeling groups.  We shall give
a presentation of $\Gamma^1$ as the fundamental group of a graph of groups.

Let $Y$ be the graph given by a semisimplicial set of the form
$$ X^{ss} \amalg \br{X}^{ss} \leftleftarrows X^{ss}_0(\ell). $$
(Here, $\br{X}^{ss}$ is $X^{ss}$ --- we have placed a bar over it to 
distinguish
the two identical factors in the coproduct.)
The face maps $d_i$ are given by
\begin{align*}
d_0 : & X^{ss}_0(\ell) \xrightarrow{t} \br{X}^{ss} \hookrightarrow X^{ss}
\amalg \br{X}^{ss},
\\
d_1 : & X^{ss}_0(\ell) \xrightarrow{s} X^{ss} \hookrightarrow X^{ss}
\amalg \br{X}^{ss}
\end{align*}
where the maps $s$ and $t$
are given on isomorphism classes by
\begin{align*}
s: & [C',H] \mapsto [C'], \\
t: & [C',H] \mapsto [C'/H]. 
\end{align*}
We give $Y$ the structure of a graph of groups $(Y,G_{(-)})$
by labeling the edges and
vertices with groups as follows.
\begin{align*}
G_{[C']} & = \aut(C') \quad \text{for $[C']$ in $X^{ss}$ or $\br{X}^{ss}$.} 
\\
G_{[C',H]} & = \aut(C',H) \quad \text{for $[C',H]$ in $X^{ss}_0(\ell)$. }
\end{align*}
We associate to the face maps $d_i$ of $Y$ monomorphisms
$$ (d_i)_* : G_{[C',H]} \rightarrow G_{d_i([C'])}. $$
The monomorphism $d_1$
is given by the natural inclusion 
$$ (d_1)_* \aut(C',H) \hookrightarrow \aut(C'). $$
Any automorphism of $C'$ which preserves a subgroup $H$ descends to an 
automorphism
of $C'/H$, and this gives the second of the two maps
$$ (d_0)_* : \aut(C',H) \rightarrow \aut(C'/H). $$

\begin{lem}
The map $(d_0)_*$ is a monomorphism.
\end{lem}

\begin{proof}
Suppose that $(d_0)_*(\alpha) = \gamma = (d_0)_*(\alpha')$ for $\alpha$ and
$\alpha'$ in $\aut(C',H)$.  Let $\phi$ be the
quotient isogeny
$C' \rightarrow C'/H$.
The automorphism $\gamma$ of $C'$ satisfies 
$$ \phi \circ \alpha = \gamma \circ \phi = \phi \circ \alpha'. $$
By composing the above equation with the dual isogeny $\widehat{\phi}$, we
see that there is an equality 
$$ \ell \cdot \alpha = \ell \cdot \alpha' $$
in the endomorphism ring $\End(C')$.  Since this ring is torsion-free, we
conclude that $\alpha = \alpha'$.
\end{proof}

The group $\Gamma^1$ acts on the tree $\calJ$ without inversions.
Proposition~\ref{prop:orbits1} shows that $Y$ is the quotient $\Gamma^1
\setminus \calJ$.
Bass-Serre theory \cite[I.5.4]{Serre}, 
combined with Proposition~\ref{prop:isotropy1} 
immediately gives the following theorem.

\begin{thm}\label{thm:amalg}
The group $\Gamma^1$ is the fundamental group of the  
graph of groups $(Y,G_{(-)})$.
\end{thm}

\section{$K(2)$-local topological modular forms}\label{sec:topmodforms}

\subsection{Morava $E$-theories}\label{sec:Etheory}

Goerss and Hopkins \cite{GoerssHopkins} refined the Hopkins-Miller Theorem
\cite{Rezk} to produce a functor
\begin{align*}
E : \mathcal{FGL} & \rightarrow \text{$E_\infty$-{\it ring spectra}} \\
(k, F) & \mapsto E(k,F).
\end{align*}
Here, $\mathcal{FGL}$ is the category of pairs $(k, F)$, where
$k$ is a perfect field of characteristic $p$ and $F$ is a formal group of 
finite height over $k$.
The spectrum $E(k,F)$ is complex orientable, and its associated formal group 
is the Lubin-Tate universal deformation $F$.

The Goerss-Hopkins-Miller functor extends naturally to the category of 
pairs $(k, F)$ obtained by 
insisting that
the ground ring $k = \prod_i k_i$ is only a product of perfect fields of
characteristic $p$, via the assignment
$$ E(k, F) = \prod_i E(k_i, F\vert_{k_i}). $$
In this paper, we are using $E_n$ to denote the spectrum $E(\br{\FF}_p,
H_n)$, where $H_n$ is the Honda height $n$ formal group.  Functoriality
gives rise to an action of the extended Morava stabilizer group
$$ \GG_n = \aut(\br{\FF}_{p^n},H_n) = \aut(H_n)\rtimes 
Gal(\br{\FF}_p/\FF_p). $$
We remark that the subgroup $Gal(\br{\FF}_p/\FF_{p^n})$ of
$Gal(\br{\FF}_p/\FF_p)$ acts trivially on $\aut(H_n)$.

Our reason for working over $\br{\FF}_{p}$ is that 
formal groups over a separably
closed field of positive characteristic $p$
are classified by their height \cite{Lazard}.  Therefore,
given $F$, a formal group of height $n$ over $\br{\FF}_p$, there is an
isomorphism $\alpha: F \cong H_n$, and hence an isomorphism of $E_\infty$-ring
spectra $E(\br{\FF}_p, F) \simeq E(\br{\FF}_p, H_n) = E_n$, which depends
on the isomorphism $\alpha$.

\subsection{Homotopy fixed points}

Because we make extensive use of homotopy fixed point constructions, we
pause to explain their meaning in the context of this paper.
Let $k$ be a finite extension of $\FF_{p^n}$.
Devinatz and Hopkins \cite{DevinatzHopkins} gave a construction of homotopy
fixed point spectra (which we shall denote $E_k^{\td{h}H}$) 
of the spectrum 
$E_k = E(k, H_n)$ 
with respect to closed 
subgroups $H$ of the 
profinite group
$$ G_k = \aut (k, H_n) = 
\aut(H_n) \rtimes Gal(k/\FF_p). $$ 
Actually, \cite{DevinatzHopkins} is written in the context of $k =
\FF_{p^n}$, but
there was nothing
in the theory of \cite{DevinatzHopkins} that prevented these authors
from replacing $\FF_{p^n}$ with the finite extension $k$. 

Goerss and Hopkins \cite{GoerssHopkins} proved that for extensions $k,k'$, 
the space of $E_\infty$-ring maps 
$E_\infty(E_k,E_{k'})$
has contractible components.  Thus the 
rectification methods of Dwyer, Kan, and Smith
\cite{DwyerKanSmith},\cite[3.2]{DevinatzHopkins} may be used to show that
the construction of the spectrum $E_k^{\td{h}H}$ may be made functorial in $k$. 

For a profinite group $G$ there is a more conventional notion of a
\emph{discrete $G$-spectrum}
that has been investigated by Thomason, Jardine, Goerss, Davis and others
(see, for instance, \cite{Davis}).  
Let $\Set_G$ be the Grothendieck site of finite discrete $G$-sets.
A discrete $G$-spectrum may be modeled as a presheaf of
spectra on this site.  The homotopy fixed points are
given by Quillen derived functors of the global sections functor with
respect to the model structure of \cite{Jardine}.
Given a closed subgroup $H$ of $G$, there is a restriction functor 
$\Res_H^G$ that takes
presheaves of spectra on $\Set_G$ to presheaves of spectra on $\Set_H$.

Following Daniel Davis \cite{Davis}, we shall regard the Devinatz-Hopkins
construction as producing a presheaf $\mathcal{E}_n(-)$ of spectra on the site 
$\Set_{\GG_n}$.
For an open subgroup $U$ of $\GG_n$, let $W(U)$ be the subgroup
$$ W(U) = 
U \cap Gal(\br{\FF}_p/\FF_{p^n}) \le Gal(\br{\FF}_p/\FF_{p^n}) \lhd \GG_n.
$$
Define $k(U)$ to be the finite fixed field
$$ k(U) = \br{\FF}_{p}^{W(U)}. $$
The value of the presheaf $\mathcal{E}(-)$ 
on the transitive finite discrete $\GG_n$-set $\GG_n/U$ 
is given by
$$ \mathcal{E}_n(\GG_n/U) = E_{k(U)}^{\td{h}U/W(U)} $$
For $H$ a closed subgroup of $\GG_n$, we 
\emph{define} the homotopy fixed point spectrum as the $K(n)$-localization
of the derived global sections of the restricted presheaf
\begin{equation}\label{eq:hfp}
E_n^{hH} = (R\Gamma \Res_H^{\GG_n} \mathcal{E}_n)_{K(n)}.
\end{equation}
Davis showed that these constructions are equivalent to those of Devinatz
and Hopkins. The statement of his 
theorem given below is a mild extension of 
the statement which appears in \cite{Davisthesis}.

\begin{thm}[Davis \cite{Davisthesis}]
There is an equivalence $E_n^{\td{h}H} \simeq E_n^{hH}$. 
\end{thm}

The Galois descent properties of $E_n$ are axiomatized by Rognes
\cite{Rognes}.  In his language, the spectrum 
$E_n$ is a $K(n)$-local profinite 
Galois extension of $S_{K(n)}$.  The homotopy fixed points of such spectra
are remarkably well behaved, as demonstrated in \cite{BehrensDavis}.  In
particular we show that, when homotopy fixed point spectra are defined in
the sense of Equation~\ref{eq:hfp}, we may iterate the homotopy fixed point
construction.

\begin{prop}[Behrens-Davis \cite{BehrensDavis}]
For $K$ a closed normal subgroup of $H$, a closed subgroup of $\GG_n$, 
there is an equivalence $(E_n^{hK})^{hH/K} \simeq E_n^{hH}$.
\end{prop}

\begin{rmk}
Devinatz has investigated a different approach to iterated homotopy fixed
points that differs philosophically from ours \cite{Devinatz}.  Namely, he
\emph{defines} the iterated fixed point construction $(E_n^{\td{h}K})^{hH/K}$ 
to be the spectrum $E_n^{\td{h}H}$ and then shows that this definition
makes sense (e.g. there is an associated Lyndon-Hochschild-Serre spectral
sequence).
\end{rmk}

Our reasons for engaging in this rhetorical yoga surrounding the
construction of homotopy fixed points is twofold.  Firstly, for $\Lambda$
a discrete group which lies as a subgroup in a profinite subgroup $G$, and
for $E$ a discrete $G$-spectrum, 
there is a natural map
$$ E^{hG} \rightarrow E^{h\Lambda} $$
where the spectrum $E^{h\Lambda}$ is the ordinary homotopy fixed
point spectrum.  Producing this map using only the Devinatz-Hopkins
language is less transparent.  Secondly, in this language, we can employ
the following lemma more freely.

\begin{thm}[Goerss {\cite[Theorem~6.1]{Goerss}}]\label{thm:Goerss}
Suppose that $E$ is a discrete $\widehat{\ZZ}$-spectrum.
Then the natural map
$$ E^{h\widehat{\ZZ}} \rightarrow E^{h\ZZ} $$
is an $H\FF_p$-equivalence.
\end{thm}

Goerss actually proved this theorem in the context of spaces, but the case
of spectra is handled by similar means, and is in some sense easier.

\subsection{Topological modular forms: an overview}\label{sec:TMF}

Let $\mathcal{M}$ be the moduli stack of generalized elliptic curves.
Goerss, Hopkins, Miller and their collaborators have constructed a sheaf 
$\mathcal{O}_{ell}$ (in
the \'etale topology) 
of $E_\infty$-ring spectra over $\mathcal{M}$.  The spectrum $\tmf$ is
given by the connective cover of the global sections
$$ \tmf = \tau_{\ge 0}\mathcal{O}_{ell}(\mathcal{M}). $$
The global sections $\mathcal{O}_{ell}(\mathcal{M})$ then give the
$E(2)$-localization $\tmf_{E(2)}$.  
Let $\mathcal{M}^{ns}$ be the substack of non-singular elliptic
curves.  The spectrum $\TMF$ is the spectrum of sections
$\mathcal{O}_{ell}(\mathcal{M}^{ns})$.  Let $\mathcal{M}^{ss}$
be a formal neighborhood of the mod $p$ supersingular locus of
$\mathcal{M}$.  Then the $K(2)$-localization $\tmf_{K(2)} = \TMF_{K(2)}$ 
is the spectrum of sections
$$ \TMF_{K(2)} = \mathcal{O}_{ell}(\mathcal{M}^{ss}). $$
Let $\mathcal{M}_0(\ell)$ be the moduli stack of
elliptic curves with $\Gamma_0(\ell)$-structures.  The forgetful
map
$$ \phi_f: \mathcal{M}_0^{ss}(\ell) \rightarrow  \mathcal{M}^{ss} $$
is \'etale, so we may evaluate $\mathcal{O}_{ell}$ on $\phi_f$ to realize
$\TMF_0(\ell)_{K(2)}$ as the spectrum of sections
$\mathcal{O}_{ell}(\mathcal{M}_0^{ss}(\ell))$.

Because a detailed account of this story does not yet exist in the
literature, we
reproduce just enough of it to give the constructions, due to Goerss, Hopkins,
Miller, and their collaborators, of $\TMF_{K(2)}$ and
$\TMF_0(\ell)_{K(2)}$ that we require.
What follows is basically a recapitulation of
a lecture of Charles Rezk on the subject in a workshop on topological modular
forms held in M\"unster, Germany in 2003.  

\subsection{The neighborhood of the supersingular locus}\label{sec:Mss}
  
Let $\WW = \WW(\br{\FF}_p)$ be the Witt ring with residue field 
$\br{\FF}_p$.
We shall first describe the stack $\mathcal{M}^{ss}
\widehat{\otimes}_{\ZZ_p} \WW$.  This formal stack is a profinite Galois
covering of the stack $\calM^{ss}$, with covering group equal to $Gal =
Gal(\br{\FF}_p/\FF_p)$.  Thus we may recover the sections of sheaves (in
the \'etale topology) over $\calM^{ss}$ from their sections over
$\calM^{ss}\widehat{\otimes}_{\ZZ_p} \WW$ by taking Galois invariants.

For each
isomorphism class $[C']$ in $X^{ss}$ we choose a representative $C'$
defined over $\br{\FF}_{p}$.
Let $k$ be the perfect ring given by the product $\prod_{X^{ss}} 
\br{\FF}_{p}$.
Let $\mathbf{C}$ be the coproduct
$$ \mathbf{C} = \coprod_{C' \in X^{ss}} C' $$
defined over $k$.  
The group $\aut(\mathbf{C}) = \prod_{C' \in X^{ss}} \aut(C')$ acts on
$\mathbf{C}$ over $k$. 
The stack $\mathcal{M}^{ss}\otimes_{\FF_p} \br{\FF}_{p}$ 
gives the supersingular points in the formal
neighborhood
$\mathcal{M}^{ss} \widehat{\otimes}_{\ZZ_p} \WW $.  Then we have
\begin{align*}
\mathcal{M}^{ss} \otimes_{\FF_p} \br{\FF}_{p} 
& = \coprod_{C' \in X^{ss}} \spec(\br{\FF}_{p})//\aut(C') \\
& = \stack \left( k, \prod_{C' \in X^{ss}} \Map(\aut(C'),\br{\FF}_p) \right).
\end{align*}
(Here $\stack(-,-)$ denotes the stackification of a Hopf algebroid 
in the \'etale topology.)  The elliptic curve $\mathbf{C}$ is the pullback
of the universal elliptic curve to the cover $\spec(k)$ of
$\mathcal{M}^{ss}$.

Let $\widehat{\mathbf{C}}$ be the formal completion of $\mathbf{C}$ at the 
identity.  Let
$\mathcal{G}$ be the Lubin-Tate universal deformation of the formal group 
$\widehat{\mathbf{C}}$ 
over $\WW(k)[[u_1]]$.  The
induced action of the group $\aut(\mathbf{C})$ on $\widehat{\mathbf{C}}$ 
extends to an
action on $\mathcal{G}$ over $\WW(k)$.  

Serre-Tate theory \cite{LubinSerreTate},
\cite{Messing} implies that the formal completion functor 
\begin{gather*}
\{ \text{deformations of $\mathbf{C}$ over $\WW(k)[[u_1]]$} \} \\
\downarrow \\
\{ \text{deformations of $\widehat{\mathbf{C}}$ over $\WW(k)[[u_1]]$} \} 
\end{gather*}
is an equivalence of categories.
Therefore, there exists a deformation $\td{\mathbf{C}}$ of $\mathbf{C}$ 
whose formal
group is the universal deformation $\mathcal{G}$.  
Lubin-Tate theory \cite{LubinTate}
implies that there are no non-trivial automorphisms of the deformation
$\mathcal{G}$ which
restrict to the identity on $\widehat{\mathbf{C}}$.
Therefore, the natural map 
$$ \aut(\td{\mathbf{C}}) \xrightarrow{\cong} \aut(\mathbf{C}) $$
is an isomorphism.

The map
$$ \chi_{\td{\mathbf{C}}}: \spf(\WW(k)[[u_1]]) \rightarrow
\mathcal{M}^{ss} \widehat{\otimes}_{\ZZ_p} \WW $$
which classifies $\td{\mathbf{C}}$ descends to a map
$$ \br{\chi}_{\td{\mathbf{C}}}: \spf(\WW(k)[[u_1]])// \aut(\td{\mathbf{C}}) 
\rightarrow \mathcal{M}^{ss} \widehat{\otimes}_{\ZZ_p} \WW $$
which is an isomorphism.  The inverse classifies
the universal deformation $\mathcal{G}$.

Now that we have a model for the formal stack 
$\calM^{ss} \widehat{\otimes}_{\ZZ_p} \WW$ defined over $\WW$, 
we may use Galois descent to recover 
the formal stack over $\ZZ_p$.  While the groupoid of
$\br{\FF}_p$-points of $\calM^{ss} \widehat{\otimes}_{\ZZ_p} \WW$
is given by supersingular elliptic curves over $\br{\FF}_p$ and isomorphisms
which cover the identity on $\br{\FF}_p$, the groupoid of
$\br{\FF}_p$-points of $\calM^{ss}$ consist of supersingular curves over
$\br{\FF}_p$ and isomorphisms which are not required to cover the identity
on $\br{\FF}_p$.

In the case of the universal supersingular elliptic curve $\mathbf{C}$, 
the extra
automorphisms arising from the Frobenius may be encoded in an action of the 
Galois group $Gal = Gal(\br{\FF}_p/\FF_p)$ on the 
groupoid $(\spec(k), \aut(\mathbf{C}))$.
Recall from Section~\ref{sec:Galois} that
for each curve $C' \in X^{ss}$, there is a Frobenius morphism
$$ \Frob_p: C' \rightarrow \sigma_*C' $$
where the curve $\sigma_*C'$ is a (possibly different) curve in $X^{ss}$.  
Thus there is a map
$$ \sigma_*: X^{ss} \rightarrow X^{ss}. $$
The action of the generator $\sigma$ on the objects $\spec(k)$ is
given by the composite
$$ \sigma^*: k = \prod_{C' \in X^{ss}} \br{\FF}_p 
\xrightarrow{\mathrm{permute}} \prod_{C'\in X^{ss}} \br{\FF}_p 
\xrightarrow{\prod \sigma} \prod_{C'\in X^{ss}} \br{\FF}_p = k.
$$
The induced map
$$ \sigma_*: \spec(k) \rightarrow \spec(k) $$
induces the action of $Gal$ on the objects of our groupoid.

The morphisms $\Frob_p$ assemble to 
give an
automorphism of $\mathbf{C}$ which covers $\sigma_*$.
$$
\xymatrix{
\mathbf{C} \ar[r]^{\Frob_p} \ar[d] &
\mathbf{C} \ar[d]
\\
\spec(k) \ar[r]_{\sigma_*} &
\spec(k)
}
$$
The action of $\sigma$ on the group $\aut(\mathbf{C})$ is given by 
conjugation by
the automorphism $\Frob_p$ of $\mathbf{C}$.  We have
$$ \sigma_*\alpha = \Frob_p \alpha \Frob_p^{-1} $$
for each $\alpha \in \aut(\mathbf{C})$.

There is a profinite Galois covering of formal stacks.
$$
\xymatrix{
\calM^{ss} \widehat{\otimes}_{\ZZ_p} \WW \ar[d]^{Gal}
\\
\calM^{ss}
}
$$
In a manner completely analogous to Section~\ref{sec:Galois},
the automorphism group $\aut(\mathbf{C})$ may be enlarged to
include the automorphism $\Frob_p$, giving rise to  
an extension
$$ \aut_{/\FF_p}(\mathbf{C}) = \aut(\mathbf{C}) \rtimes Gal. $$

\subsection{Construction of $\TMF_{K(2)}$}\label{sec:TMF_K(2)}

As described in Section~\ref{sec:Etheory}, 
the Goerss-Hopkins-Miller Theorem
gives an $E_\infty$-ring spectrum
$$ E(k,\widehat{\mathbf{C}}) \cong \prod_{C' \in X^{ss}} E(\br{\FF}_p, C') $$
and an action of the group $\aut_{/\FF_p}(\mathbf{C})$ on this spectrum by 
$E_\infty$-ring maps.  The coefficient ring of this complex orientable spectrum 
is given by
$$ E(k,\widehat{\mathbf{C}})_* = \WW(k)[[u_1]][u^{\pm 1}] $$
where $\abs{u} = -2$.  

The spectrum $\TMF_{K(2)}$ is defined to be the homotopy fixed
point spectrum
\begin{align*}
\TMF_{K(2)} 
& = E(k,\widehat{\mathbf{C}})^{h\aut_{/\FF_p}(\mathbf{C})} \\
& = 
\left( \prod_{C' \in X^{ss}} E(\br{\FF}_p,\widehat{C}')^{h\aut(C')} 
\right)^{hGal}.
\end{align*}
In Section~\ref{sec:J2ell} we shall find it useful to work with a version
of $\TMF_{K(2)}$ where we do not take Galois fixed points.  We thus make the
definition
\begin{align*} 
\TMF_{K(2), \br{\FF}_p} & = E(k,\widehat{\mathbf{C}})^{h\aut(\mathbf{C})} \\
& = \prod_{C' \in X^{ss}} E(\br{\FF}_p,\widehat{C}')^{h\aut(C')}.
\end{align*}

\subsection{$\Gamma_0(\ell)$-structures}\label{sec:TMF_0(ell)}

The construction of $\TMF_0(\ell)$ is completely analogous.  One simply
replaces everywhere the formal moduli stack $\mathcal{M}^{ss}$ with
$\mathcal{M}_0(\ell)^{ss}$.  The automorphism groups $\aut(C')$ are replaced
with $\aut(C',H)$ for $(C',H) \in X_0^{ss}(\ell)$.

Explicitly, let $k'$ be the perfect ring
$$
k' 
= \prod_{\Gamma_0(\ell)(\mathbf{C})} k
= \prod_{X_0^{ss}(\ell)} \br{\FF}_{p}
$$
where $\Gamma_0(\ell)(\mathbf{C})$ is the set of $\Gamma_0(\ell)$-structures on
$\mathbf{C}$.  
We define $\mathbf{C}_0(\ell)$ to be the elliptic curve over $k'$ given by
$$ \mathbf{C}_0(\ell) = \coprod_{\Gamma_0(\ell)(\mathbf{C})} \mathbf{C}. $$
We give $\mathbf{C}_0(\ell)$ the canonical $\Gamma_0(\ell)$-structure
$\mathbf{H}$ which restricts to $H$ on the component 
corresponding to 
the element 
$H \in \Gamma_0(\ell)(\mathbf{C}) $.

Since the map
$$ \mathbf{C}[\ell] \rightarrow \spec(k) $$
is \'etale, given a $\Gamma_0(\ell)$-structure $H$ on $\mathbf{C}$, 
there is a unique
extension to a $\Gamma_0(\ell)$-structure $\td{H}$ on $\td{\mathbf{C}}$ over
$\WW(k)[[u_1]]$.
The elliptic curve over $\WW(k')[[u_1]]$ given by
$$ \td{\mathbf{C}}_0(\ell) = \coprod_{\Gamma_0(\ell)(\mathbf{C})} \td{C} $$
is a deformation of $\mathbf{C}_0(\ell)$.
The $\Gamma_0(\ell)$-structure $\mathbf{H}$ extends uniquely to
a $\Gamma_0(\ell)$-structure $\td{\mathbf{H}}$ on 
$\td{\mathbf{C}}_0(\ell)$. 
It restricts to $\td{H}$ on the component corresponding to 
the element $H \in \Gamma_0(\ell)(\mathbf{C})$.

Define the group $\aut(\mathbf{C}_0(\ell), \mathbf{H})$ to be the
finite group of
automorphisms of $\mathbf{C}_0(\ell)$ which preserve the level structure
$\mathbf{H}$:
$$ \aut(\mathbf{C}_0(\ell), \mathbf{H}) = \prod_{(C',H) \in
X^{ss}_0(\ell)} \aut(C', H). $$
The automorphism $\Frob_p$ on $\mathbf{C}$ of Section~\ref{sec:Mss} will 
permute
the $\Gamma_0(\ell)$-structures, inducing an action of $Gal$ on the groupoid
$(\spec(k'), \aut(\mathbf{C}_0(\ell), \mathbf{H}))$.  We get an
extension of groups
$$ \aut_{/\FF_p}(\mathbf{C}_0(\ell), \mathbf{H}) = \aut(\mathbf{C}_0(\ell),
\mathbf{H})\rtimes Gal.$$

Just as in the case of $\TMF_{K(2)}$, we use the Goerss-Hopkins-Miller 
theorem to
produce a spectrum $E(k',\widehat{\mathbf{C}}_0(\ell))$.
The spectrum $\TMF_0(\ell)_{K(2)}$ 
is given as follows:
\begin{align*}
\TMF_0(\ell)_{K(2)}
& = E(k',\widehat{\mathbf{C}}_0(\ell))^{h\aut_{/\FF_p}(\mathbf{C}_0(\ell),
\mathbf{H})} \\
& = \left( \prod_{(C',H) \in X^{ss}_0(\ell)} E(\br{\FF}_p,
\widehat{C}')^{h\aut(C',H)} \right)^{hGal}.
\end{align*}

Just as in the case of $\TMF$, we will define $\TMF_0(\ell)_{K(2),\br{\FF}_p}$ 
to be the version where we do not take Galois fixed points:
\begin{align*}
\TMF_0(\ell)_{K(2),\br{\FF}_p}
& = E(k',\widehat{\mathbf{C}}_0(\ell))^{h\aut(\mathbf{C}_0(\ell),
\mathbf{H})} \\
& = \prod_{(C',H) \in X^{ss}_0(\ell)} E(\br{\FF}_p,
\widehat{C}')^{h\aut(C',H)}.
\end{align*}

\section{Relation to the spectrum $Q(\ell)$}\label{sec:Qell}

\subsection{The spectrum $Q(\ell)$}\label{sec:Q(ell)const}

In \cite{BehrensK(2)}, using the sheaf $\mathcal{O}_{ell}$ of
Section~\ref{sec:TMF}, we introduced a spectrum $Q(\ell)$ built out of 
$\TMF$ and $\TMF_0(\ell)$.  We give an independent $K(2)$-local construction 
here.  Nevertheless, the reader might find it useful to
refer to \cite{BehrensK(2)},
where the motivation for the construction is given.

The spectrum $Q(\ell)_{K(2)}$ is the totalization of a 
semi-cosimplicial $E_\infty$-ring spectrum
of the following form:
\begin{equation}\label{eq:cosimpQ(ell)}
\TMF_{K(2)} \quad \Rightarrow \quad
\begin{array}{c}
\TMF_0(\ell)_{K(2)} \\
\times \\
\TMF_{K(2)}
\end{array}
\quad \Rrightarrow \quad
\TMF_0(\ell)_{K(2)}.
\end{equation}
The coface maps are given in terms of certain maps of $E_\infty$-ring
spectra:
\begin{align*}
\phi_f^* : & \TMF_{K(2)} \rightarrow \TMF_0(\ell)_{K(2)}, \\
\phi_q^* : & \TMF_{K(2)} \rightarrow \TMF_0(\ell)_{K(2)}, \\
\psi_{[\ell]}^* : & \TMF_{K(2)} \rightarrow \TMF_{K(2)}, \\
\psi_d^* : & \TMF_0(\ell)_{K(2)} \rightarrow \TMF_0(\ell)_{K(2)}.
\end{align*}
The coface maps on $0$-cosimplicies 
$$ d_i: \TMF_{K(2)} \rightarrow \TMF_0(\ell)_{K(2)} \times \TMF_{K(2)} $$ 
are defined by
\begin{gather*}
d_0 = \phi_q^* \times \psi_{[\ell]}^*, \\
d_1 = \phi_f^* \times Id. 
\end{gather*}
The coface maps 
$$ d_i: \TMF_0(\ell)_{K(2)} \times 
\TMF_{K(2)} \rightarrow \TMF_0(\ell)_{K(2)} $$ 
are defined by
\begin{gather*}
d_0 = \psi_d^* \circ p_1, \\
d_1 = \phi_f^* \circ p_2, \\
d_2 = p_1
\end{gather*}
where $p_1$, $p_2$ are the projections onto the first and second factors of
the product $\TMF_0(\ell)_{K(2)} \times \TMF_{K(2)}$.

We produce the required maps using the Goerss-Hopkins-Miller 
functor (Section~\ref{sec:Etheory}).
\vspace{10pt}

\noindent
{\bf The map $\psi_{[\ell]}^*$:}  
The $\ell$th power isogeny
$$ [\ell] : \mathbf{C} \rightarrow \mathbf{C} $$
induces an automorphism
$$ \psi_{[\ell]} = (Id,[\ell]): (k,\widehat{\mathbf{C}}) \rightarrow
(k,\widehat{\mathbf{C}}). $$
Applying the Goerss-Hopkins-Miller functor, we get a map
$$ \psi_{[\ell]}^*: E(k,\widehat{\mathbf{C}}) \rightarrow
E(k,\widehat{\mathbf{C}}). $$
Because the $\ell$-power isogeny commutes with all of the automorphisms of
$\mathbf{C}$, this map descends to the homotopy fixed points
\begin{eqnarray*}
\psi_{[\ell]}^*: \TMF_{K(2)}
& = & E(k,\widehat{\mathbf{C}})^{h\aut_{/\FF_p}(\mathbf{C})} \\
& \rightarrow & E(k,\widehat{\mathbf{C}})^{h\aut_{/\FF_p}(\mathbf{C})} \\
& = & \TMF_{K(2)}.
\end{eqnarray*}
We remark that by replacing the pair $(k,\mathbf{C})$ with the pair
$(k',\mathbf{C}_0(\ell))$, we get a map
$$ \psi_{[\ell]}^*: \TMF_0(\ell)_{K(2)} 
\rightarrow \TMF_0(\ell)_{K(2)}.$$

\noindent
{\bf The map $\psi_d^*$:}
In Section~\ref{sec:J'face} we defined,
for each pair $(C',H) \in X^{ss}_0(\ell)$, 
a pair $(C_H,d(H)) \in X^{ss}_0(\ell)$, and a
degree $\ell$ 
isogeny
$$ \phi_H : C' \rightarrow C_H. $$
The isogeny $\phi_H$ has kernel $H$, and the dual isogeny
$\widehat{\phi}_H$ has kernel $d(H)$.
Observe that the pair $(C_H,d(H))$ actually determines $(C',H)$: we have
\begin{align*}
C_{d(H)} & = C', \\
d(d(H)) & = H.
\end{align*}

We may define an involution
$$ \br{\psi}_d : k' = \prod_{(C',H) \in X_0^{ss}(\ell)} \br{\FF}_{p} 
\rightarrow 
\prod_{(C',H) \in X_0^{ss}(\ell)} \br{\FF}_{p} = k'
$$
given by permuting the factors: we send the factor corresponding to 
$(C',H)$ to the
factor corresponding to $(C_H, d(H))$.
The maps $\phi_H$ assemble to give a degree $\ell$ isogeny
$$ \td{\psi}_d: \mathbf{C}_0(\ell) \rightarrow \mathbf{C}_0(\ell) $$
which covers the map $\br{\psi}_d$.
The kernel of $\td{\psi}_d$ is $\mathbf{H}$.

Since
$\phi_{d(H)}$ and the dual isogeny $\widehat{\phi}_H$ have the same
kernel $d(H)$, there exists an automorphism $\gamma_H$ of $C'$ over
$\br{\FF}_{p}$ so that the
following diagram commutes.
\begin{equation}\label{diag:relphiH}
\xymatrix{
C_H \ar[r]^{\widehat{\phi}_H} \ar[dr]_{\phi_{d(H)}} &
C' \ar@{.>}[d]^{\gamma_H} 
\\
& C'
}
\end{equation}
By applying the dual isogeny functor to the above diagram, we see that
$\gamma_H$ preserves $H$, so $\gamma_H$ actually lies in the automorphism group
$\aut(C',H)$.
Diagram~(\ref{diag:relphiH}) gives us the relation
\begin{equation}\label{eq:relphiH}
\phi_{d(H)} \circ \phi_H = [\ell] \circ \gamma_H.
\end{equation}
The automorphisms $\gamma_H$ assemble to give an automorphism
$$ \gamma: \mathbf{C}_0(\ell) \rightarrow \mathbf{C}_0(\ell)
$$
defined over $k'$
which preserves the subgroup $\mathbf{H}$.
Equation~(\ref{eq:relphiH}) gives us the relation
\begin{equation}\label{eq:relpsid}
\td{\psi}_d \circ \td{\psi}_d = [\ell] \circ \gamma.
\end{equation}

We assemble these automorphisms to get an automorphism of pairs
$$ \psi_d = (\br{\psi}_d, (\td{\psi}_d)_*): 
(k',\widehat{\mathbf{C}}_0(\ell)) 
\rightarrow (k',\widehat{\mathbf{C}}_0(\ell)) $$
which induces a map
$$ \psi_d^*: E(k',\widehat{\mathbf{C}}_0(\ell)) 
\rightarrow E(k',\widehat{\mathbf{C}}_0(\ell)). $$
We claim that $\psi_d^*$ descends to an automorphism of $\TMF_0(\ell)_{K(2)}$.
Suppose that $\beta$ is an element of $\aut_{/\FF_p}(\mathbf{C}_0(\ell),
\mathbf{H})$.  Then we have
the following diagram.
\begin{equation}\label{diag:betabar}
\xymatrix{
\mathbf{C}_0(\ell) \ar[r]^{\td{\psi}_d} \ar[d]_{\beta} &
\mathbf{C}_0(\ell)
\ar@{.>}[d]^{\br{\beta}}
\\
\mathbf{C}_0(\ell) \ar[r]_{\td{\psi}_d} &
\mathbf{C}_0(\ell)
}
\end{equation}
Since $\beta$ preserves $\mathbf{H} = \ker \td{\psi}_d $, 
it descends uniquely
to give $\br{\beta}$. Applying the dual isogeny functor to
Diagram~(\ref{diag:betabar}), we see that the map $\br{\beta}$ preserves the
kernel of the dual isogeny of $\td{\psi}_d$, But we have argued that this
kernel is 
also given by $\mathbf{H}$.  Thus
$\br{\beta}$ also lies in $\aut_{/\FF_p}(\mathbf{C}_0(\ell), \mathbf{H})$.  
We conclude that
the map $\psi_d^*$ descends to the 
$\aut_{/\FF_p}(\mathbf{C}_0(\ell), \mathbf{H})$-fixed points
$ \TMF_0(\ell)_{K(2)} = 
E(k',\widehat{\mathbf{C}}_0(\ell))^{h\aut_{/\FF_p}(\mathbf{C}_0(\ell), 
\mathbf{H})}$ to give a map
$$ \psi_d^*: 
\TMF_0(\ell)_{K(2)} \rightarrow 
\TMF_0(\ell)_{K(2)}. $$
Because $\gamma$ is contained in 
$\aut_{/\FF_p}(\mathbf{C}_0(\ell), \mathbf{H})$, 
it acts trivially on $\TMF_0(\ell)_{K(2)}$, and we have
the following relation on $\TMF_0(\ell)_{K(2)}$.
\begin{equation}\label{eq:psi_drel}
\psi_d^* \circ \psi_d^* = \psi_{[\ell]}^*. 
\end{equation}

\begin{rmk}
The equality in Equation~(\ref{eq:psi_drel}) is a strict equality that
occurs on the point-set level.  This is because the
homotopy fixed point spectrum is the \emph{actual} fixed points of an 
appropriate
fibrant replacement.
\end{rmk}

\noindent
{\bf The map $\phi_f^*$:}
Let $\br{\chi}$ denote the diagonal map
$$ \br{\chi} : k \rightarrow \prod_{\Gamma_0(\ell)(\mathbf{C})} k = k'. $$
Over this map we have $\mathbf{C}_0(\ell) = \mathbf{C} \otimes_k k'$.
We therefore get a map of pairs
$$ \chi: (k',\widehat{\mathbf{C}}_0(\ell)) \rightarrow
(k,\widehat{\mathbf{C}}).
$$
The diagonal embedding 
$$
\aut_{/\FF_p}(\mathbf{C}) \hookrightarrow \aut_{/\FF_p}(\mathbf{C}_0(\ell)) 
$$
is compatible with the map
$\chi$.
The natural inclusion 
$$ 
\iota: \aut_{/\FF_p}(\mathbf{C}_0(\ell),\mathbf{H})
\hookrightarrow \aut_{/\FF_p}(\mathbf{C}_0(\ell))
$$
gives $\phi_f^*$ as the composite
\begin{align*}
\phi_f^*: \TMF_{K(2)} 
& = E(k,\widehat{\mathbf{C}})^{h\aut_{/\FF_p}(\mathbf{C})} \\
& \xrightarrow{\chi^*} E(k',\widehat{\mathbf{C}}_0(\ell))^{h\aut_{/\FF_p}(
\mathbf{C}_0(\ell))} \\
& \xrightarrow{\iota^*} 
E(k',\widehat{\mathbf{C}}_0(\ell))^{h\aut_{/\FF_p}(\mathbf{C}_0(\ell),
\mathbf{H})} \\
& = \TMF_0(\ell)_{K(2)}.
\end{align*}
The commutativity of the diagram
$$
\xymatrix{
(k',\widehat{\mathbf{C}}_0(\ell)) \ar[r]^\chi \ar[d]_{\psi_{[\ell]}} &
(k,\widehat{\mathbf{C}}) \ar[d]^{\psi_{[\ell]}}
\\
(k',\widehat{\mathbf{C}}_0(\ell)) \ar[r]_\chi &
(k,\widehat{\mathbf{C}})
} $$
implies the relation 
\begin{equation}\label{eq:psi_ellrel}
\phi_f^* \psi_{[\ell]}^* = \psi_{[\ell]}^* \phi_f^*.
\end{equation}

\noindent
{\bf The map $\phi_q^*$:}
The map $\phi_q^*$ is defined to be the composite
$$ \phi_q^*: \TMF_{K(2)} \xrightarrow{\phi_f^*} \TMF_0(\ell)_{K(2)}
\xrightarrow{\psi_d^*} \TMF_0(\ell)_{K(2)}. $$

The construction of the spectrum $Q(\ell)_{K(2)}$ 
is completed by the following
lemma.

\begin{lem}
The coface maps in (\ref{eq:cosimpQ(ell)}) satisfy the cosimplicial
identities.
\end{lem}

\begin{proof}
We translate the cosimplicial identities into the maps that define the
$d_i$'s.
\begin{alignat}{2}
d_0 d_0 & = d_1 d_0 & \qquad 
\psi_d^* \phi_q^* & = \phi_f^* \psi_{[\ell]}^*, 
\label{eq:cosimp1}
\\
d_2 d_0 & = d_0 d_1 & \qquad
\phi_q^* & = \psi_d^* \phi_f^*, 
\label{eq:cosimp2}
\\
d_1 d_1 & = d_2 d_1 & \qquad
\phi_f^* & = \phi_f^*.
\label{eq:cosimp3}
\end{alignat}
Relation~(\ref{eq:cosimp3}) is tautologous, and Relation~(\ref{eq:cosimp2})
is immediate from our definition of $\phi_q^*$.
Relation~(\ref{eq:cosimp1}) then follows from Relation~(\ref{eq:cosimp2}),
Equation~(\ref{eq:psi_drel}), and Equation~(\ref{eq:psi_ellrel}).
\end{proof}

\subsection{$Q(\ell)_{K(2)}$ as the homotopy fixed point spectrum
$E(\Gamma)$}\label{sec:J2ell}

In this section we shall prove the following theorem.

\begin{thm}\label{thm:Q(ell)}
There is an equivalence $Q(\ell) \simeq E(\Gamma) = (E_2^{h\Gamma_{Gal}})$.
\end{thm}

Before we prove these theorems we address some finer points concerning
our use of 
Morava $E$-theory.  We have fixed a supersingular elliptic 
curve $C$ defined over
$\br{\FF}_p$, and have fixed an isomorphism between it and the Honda height 
$2$ formal group $H_2$ over $\br{\FF}_p$.  This gives rise, by the
Goerss-Hopkins-Miller theorem, to a fixed
isomorphism
$$ E(\br{\FF}_p, \widehat{C}) \cong E(\br{\FF}_p, H_2) = E_2. $$
Our fixed isomorphism $\widehat{C} \cong H_2$ also gives an isomorphism
$$ \aut(\br{\FF}_p, \widehat{C}) \cong \aut(\br{\FF}_p, H_2) = \GG_2. $$
In what follows, when we refer to $E_2$ and $\GG_2$, we shall actually be
implicitly \emph{identifying} these with 
$E(\br{\FF}_p,\widehat{C})$ and $\aut(\br{\FF}_p,\widehat{C})$ using our
fixed isomorphisms.

We recall how the Goerss-Hopkins-Miller theorem gives rise to an action of
$\GG_2$ on $E_2$ by $E_\infty$-ring maps.  Let $g$ be an element of $\GG_2$.  It
is an automorphism
$$ g = (g_0,g_1): (\br{\FF}_p,\widehat{C}) \rightarrow
(\br{\FF}_p,\widehat{C}). $$
Because the Goerss-Hopkins-Miller theorem gives a \emph{contravariant}
functor $E(-,-)$, the left action $L_g$ of $g$ on $E_2$ is given by the 
image of 
$g^{-1}$
under the functor $E(-,-)$
\begin{equation}\label{eq:L_g}
L_g = (g^{-1})^* : E_2 \rightarrow E_2.
\end{equation}

For $C' \in X^{ss}$, let $E_{C'}$ denote 
the spectrum $E(\br{\FF}_p,\widehat{C}')$.
We defined $\TMF_{K(2)}$ and $\TMF_0(\ell)_{K(2)}$ as homotopy fixed points
of the spectra
\begin{align*}
E(k, \widehat{\mathbf{C}}) & \cong \prod_{C' \in X^{ss}} E_{C'},
\\
E(k', \widehat{\mathbf{C}}_0(\ell)) & \cong \prod_{(C',H) \in
X_0^{ss}(\ell)} E_{C'}.
\end{align*}
The spectra $E_{C'}$ are isomorphic to $E_2 =
E_{C}$ using the fixed isomorphisms (over $\br{\FF}_p$)
of formal groups 
$$ (\phi_{C'})_* : \widehat{C} \rightarrow \widehat{C}' $$
induced by the isogenies $\phi_{C'}$ of Section~\ref{sec:orbits}.  We get
an induced isomorphism
$$ \phi_{C'}^* : E_{C'} \xrightarrow{\cong} E_2. $$
Under this isomorphism, the induced action of the group $\aut(C')$ on $E_2$
corresponds to the action given by the 
embedding $\iota_{C'}$ of $\aut(C')$ in $\Gamma$ defined
in Section~\ref{sec:orbits}.

Let $Q(\ell)_{K(2),\br{\FF}_p}$ be the spectrum obtained by the totalization
of the Galois equivariant semi-cosimplicial spectrum
$$ \TMF_{K(2),\br{\FF}_p} \Rightarrow \TMF_{K(2),\br{\FF}_p} \times
\TMF_0(\ell)_{K(2),\br{\FF}_p} \Rrightarrow \TMF_0(\ell)_{K(2),\br{\FF}_p} $$
where we have \emph{not} taken Galois fixed points.  

Let $\sigma^\ZZ \subset Gal$ be the discrete group given by powers of the
Frobenius $\sigma$.
The following lemma is
a consequence of Theorem~\ref{thm:Goerss}.

\begin{lem}\label{lem:Goerss}
The natural maps
\begin{align*}
\TMF_{K(2)} & = (\TMF_{K(2),\br{\FF}_p})^{hGal} \rightarrow
(\TMF_{K(2),\br{\FF}_p})^{h\sigma^\ZZ}, \\
\TMF_0(\ell)_{K(2)} & = (\TMF_0(\ell)_{K(2),\br{\FF}_p})^{hGal} \rightarrow
(\TMF_0(\ell)_{K(2),\br{\FF}_p})^{h\sigma^\ZZ}
\end{align*}
are equivalences.
\end{lem}

\begin{cor}
There is an equivalence
$$ Q(\ell)_{K(2)} \rightarrow (Q(\ell)_{K(2),\br{\FF}_p})^{h\sigma^\ZZ}. $$
\end{cor}

The remainder of this section is devoted to proving
Theorem~\ref{thm:Q(ell)}.
We shall first prove
that there is an equivalence
\begin{equation}\label{eq:Q(ell)Fpbar}
Q(\ell)_{K(2),\br{\FF}_p} \xrightarrow{\simeq} E_2^{h\Gamma}.
\end{equation}
We will then prove that this equivalence commutes with the action of the
Frobenius $\sigma$.  Theorem~\ref{thm:Q(ell)} is then recovered by applying
the functor $(-)^{h\sigma^\ZZ}$ to (\ref{eq:Q(ell)Fpbar}).

\begin{prop}\label{prop:Ehgammacosimp}
The homotopy fixed point spectrum $E_2^{h\Gamma}$ is the totalization of a
semi-cosimplicial spectrum of the form
\begin{equation}\label{eq:Ehgammacosimp}
\prod_{C' \in X^{ss}} E_2^{h\aut(C')} \Rightarrow
\begin{array}{c}
\prod_{(C',H) \in X^{ss}_0(\ell)} E_2^{h\aut(C',H)} \\
\times \\
\prod_{C' \in X^{ss}} E_2^{h\aut(C')}
\end{array}
\Rrightarrow 
\prod_{(C',H) \in X^{ss}_0(\ell)} E_2^{h\aut(C',H)}.
\end{equation}
The coface maps
$$ d_i : \prod_{C' \in X^{ss}} E_2^{h\aut(C')} \rightarrow
\prod_{(C',H) \in X^{ss}_0(\ell)} E_2^{h\aut(C',H)} \quad
\times \quad
\prod_{C' \in X^{ss}} E_2^{h\aut(C')} $$
are defined on components by
\begin{alignat*}{2}
(d_0)_{(C',H)} & = L_{g_{(C',H)}} \circ \Res_{\aut(C_H,d(H))}^{\aut(C_H)}
& : \qquad & E_2^{h\aut(C_H)} \rightarrow E_2^{h\aut(C',H)}, 
\\
(d_1)_{(C',H)} & = \Res_{\aut(C',H)}^{\aut(C')} 
& : \qquad & E_2^{h\aut(C')} \rightarrow E_2^{h\aut(C',H)}, 
\\
(d_0)_{C'} & = L_{\ell^{-1}}
& : \qquad & E_2^{h\aut(C')} \rightarrow E_2^{h\aut(C')}, 
\\
(d_1)_{C'} & = Id
& : \qquad & E_2^{h\aut(C')} \rightarrow E_2^{h\aut(C')}.
\end{alignat*}
The coface maps 
$$
d_i: 
\prod_{(C',H) \in X^{ss}_0(\ell)} E_2^{h\aut(C',H)} \quad \times 
\quad \prod_{C' \in X^{ss}} E_2^{h\aut(C')}
\rightarrow 
\prod_{(C',H) \in X^{ss}_0(\ell)} E_2^{h\aut(C',H)}
$$
are defined on components by
\begin{alignat*}{2}
(d_0)_{(C',H)} & = L_{g_{(C',H)}}
& : \qquad & E_2^{h\aut(C_H,d(H))} \rightarrow E_2^{h\aut(C',H)},
\\
(d_1)_{(C',H)} & = \Res_{\aut(C',H)}^{\aut(C')}
& : \qquad & E_2^{h\aut(C')} \rightarrow E_2^{h\aut(C',H)},
\\
(d_2)_{(C',H)} & = Id
& : \qquad & E_2^{h\aut(C',H)} \rightarrow E_2^{h\aut(C',H)}.
\end{alignat*}
Here the element $g_{(C',H)}$ is the element of $\Gamma$ 
defined by Equation~(\ref{eq:g_C'H}).
\end{prop}

\begin{proof}
Since the complex $\calJ'$ is non-equivariantly contractible, the natural map
$$ E_2^{h\Gamma} \rightarrow \Map(\calJ',E_2)^{h\Gamma} $$
is an equivalence.
Using the semi-simplicial structure of $\calJ'$, we see that $E^{h\Gamma}$
is equivalent to the totalization of the semi-cosimplicial spectrum
\begin{equation}\label{eq:calJcosimp} 
\Map(\calJ'_0, E_2)^{h\Gamma} \Rightarrow \Map(\calJ'_1,E_2)^{h\Gamma} 
\Rrightarrow
\Map(\calJ'_2,E_2)^{h\Gamma}.
\end{equation}
By ``Shapiro's lemma'', for a subgroup $F$ of $\Gamma$, there is an
equivalence
\begin{equation}\label{eq:Shapiro}
E_2^{hF} \simeq \Map(\Gamma/F, E_2)^{h\Gamma}.
\end{equation}
The semi-cosimplicial spectrum given in the proposition 
is obtained by substituting 
the descriptions of $\calJ'_i$ given
in Equation~(\ref{eq:J'gamma}) into
Equation~(\ref{eq:calJcosimp}), and then applying Equation~(\ref{eq:Shapiro}).  
The descriptions of the coface maps given
in the proposition follow immediately from the descriptions of the face maps of
$\calJ'_\bullet$ given by Proposition~\ref{prop:J'face}. 
\end{proof}

\noindent
{\it Construction of the equivalence (\ref{eq:Q(ell)Fpbar}).}
We first describe the maps $\psi_{[\ell]}^*$, $\psi_d^*$, $\phi_f^*$, and
$\phi_q^*$ of Section~\ref{sec:Q(ell)const} under the isomorphisms
\begin{align*}
\TMF_{K(2),\br{\FF}_p} & \cong \prod_{C' \in X^{ss}}
E_{C'}^{h\aut(C')}, \\
\TMF_0(\ell)_{K(2),\br{\FF}_p} & \cong \prod_{(C',H) \in X_0^{ss}(\ell)}
E_{C'}^{h\aut(C',H)}.
\end{align*}
We describe the components of our maps below, which are read off from their
definitions in Section~\ref{sec:Q(ell)const}:
\begin{alignat*}{2}
(\psi_{[\ell]}^*)_{C'} & = [\ell]^* 
& : \qquad & E_{C'}^{h\aut(C')} \rightarrow
E_{C'}^{h\aut(C')},
\\
(\psi_d^*)_{(C',H)} & = \phi_H^*
& : \qquad & E_{C_H}^{h\aut(C_H,d(H))} \rightarrow 
E_{C'}^{h\aut(C',H)},
\\
(\phi_f^*)_{(C',H)} & = \Res_{\aut(C',H)}^{\aut(C')}
& : \qquad & E_{C'}^{h\aut(C')} \rightarrow 
E_{C'}^{h\aut(C',H)},
\\
(\phi_q^*)_{(C',H)} & = \phi_H^* \circ \Res_{\aut(C_H,d(H))}^{\aut(C_H)}
& : \qquad & E_{C_H}^{h\aut(C_H)} \rightarrow
E_{C'}^{h\aut(C',H)}.
\end{alignat*}
The left $\MS_2$ action on $E_2$ given by Equation~(\ref{eq:L_g}) gives the
following commutative diagrams.
\begin{equation}\label{diag:essential}
\xymatrix{
E_{C'} \ar[r]^{[\ell]^*} \ar[d]_{\phi_{C'}^*} &
E_{C'} \ar[d]^{\phi_{C'}^*} &
E_{C_H} \ar[r]^{\phi_H^*} \ar[d]_{\phi_{C_H}^*} &
E_{C'} \ar[d]^{\phi_{C'}^*}
\\
E_2 \ar[r]_{L_{\ell^{-1}}} &
E_2 &
E_2 \ar[r]_{L_{g_{(C',H)}}} &
E_2
}
\end{equation}
These diagrams, using Proposition~\ref{prop:Ehgammacosimp}, give rise the
following
equivalence of semi-cosimplicial spectra.
\begin{equation}\label{diag:topbottom}
\xymatrix@R+1em{
\prod_{X^{ss}} E_{C'}^{h\aut(C')} 
\ar@<.5ex>[r] \ar@<-.5ex>[r] \ar[d]_{(\phi_{C'}^*)} & 
{\begin{array}{c}
\prod_{X^{ss}_0(\ell)} E_{C'}^{h\aut(C',H)} \\
\times \\
\prod_{X^{ss}} E_{C'}^{h\aut(C')}
\end{array}}
\ar@<1ex>[r] \ar[r] \ar@<-1ex>[r] \ar[d]|{(\phi_{C'}^*) \times (\phi_{C'}^*)} 
& 
\prod_{X^{ss}_0(\ell)} E_{C'}^{h\aut(C',H)} \ar[d]^{(\phi_{C'}^*)}
\\
\prod_{X^{ss}} E_2^{h\aut(C')} 
\ar@<.5ex>[r] \ar@<-.5ex>[r] &
{\begin{array}{c}
\prod_{X^{ss}_0(\ell)} E_2^{h\aut(C',H)} \\
\times \\
\prod_{X^{ss}} E_2^{h\aut(C')}
\end{array}}
\ar@<1ex>[r] \ar[r] \ar@<-1ex>[r] &
\prod_{X^{ss}_0(\ell)} E_2^{h\aut(C',H)}
}
\end{equation}
The totalization of the top row gives $Q(\ell)_{K(2),\br{\FF}_p}$, while 
Proposition~\ref{prop:Ehgammacosimp} implies that the
totalization of the bottom row gives $E_2^{h\Gamma}$.

We will finish this section by proving the
following lemma.

\begin{lem}\label{lem:verticalGal}
The maps $(\phi_{C'}^*)$ of Diagram~\ref{diag:topbottom} are Galois equivariant.
\end{lem}

We pause to explain how Lemma~\ref{lem:verticalGal} completes the proof of
Theorem~\ref{thm:Q(ell)}.
The coface maps of the 
top row of
Diagram~\ref{diag:topbottom} are Galois equivariant by construction.
Up to this point, we have not addressed the Galois equivariance of
the coface maps of the 
bottom row of Diagram~\ref{diag:topbottom}.
By functoriality, the vertical maps
$(\phi_{C'}^*)$ are \emph{isomorphisms} of spectra, so the Galois equivariance 
of the coface maps of the bottom row will follow from
Lemma~\ref{lem:verticalGal}.  
The proof of Theorem~\ref{thm:Q(ell)} is then completed by applying Galois
fixed points to the equivalence (\ref{eq:Q(ell)Fpbar}).

\begin{rmk}
One could have also deduced the Galois
equivariance of the coface maps of the 
bottom row from the fact that the Galois action 
on the Tate module $V_\ell(C)$ induces a Galois action on the building $\calJ'$.
\end{rmk}

\begin{proof}[Proof of Lemma~\ref{lem:verticalGal}]
We will only prove that the map 
\begin{equation}\label{eq:phiC'}
(\phi_{C'}^*) : \prod_{C' \in X^{ss}}
E_{C'}^{h\aut(C')} \rightarrow
\prod_{C' \in X^{ss}} E_2^{h\aut(C')}
\end{equation}
of the first column is
Galois equivariant.  The other case, with level structure, proceeds in the
same manner.  

We first recall the Galois action on the source and target in 
Equation~(\ref{eq:phiC'}).  By Equation~{(\ref{eq:L_g})}, the Frobenius 
$\sigma$ acts on the source by the map
induced by $(\Frob_p^{-1})^*$ (Section~\ref{sec:Galois}).
$$ \sigma : E_{C'}^{h\aut(C')} 
\xrightarrow{(\Frob^{-1}_p)^*} E_{\sigma_*C'}^{h\aut(\sigma_*C')}. $$
The Frobenius action on the target in Equation~(\ref{eq:phiC'}) is induced from the
Frobenius action on the $\Gamma$-set $\coprod_{X^{ss}}
\Gamma/\phi_{C'}^{-1}\aut(C')\phi_{C'}$ (as described by
Equation~(\ref{eq:GalJ'}))
through our application of Shapiro's lemma (Equation~(\ref{eq:Shapiro})).
The 
resulting Frobenius action is given by
$$ \sigma: E_2^{h\aut(C')} \xrightarrow{y_{C'}^* \circ (\Frob_p^{-1})^*}
E_2^{h\aut(\sigma_*C')}
$$
where the quasi-isogeny $y_{C'}$ is defined in the proof of
Theorem~\ref{thm:calJ'}.  The lemma now follows from the commutativity of
the following diagram, which is immediate given the definition of $y_{C'}$.
$$
\xymatrix@C+4em{
E_{C'}^{h\aut(C')} \ar[r]^{(\Frob_p^{-1})^*} \ar[d]_{\phi_{C'}^*} 
& E_{\sigma_*C'}^{h\aut(\sigma_*C')} \ar[d]^{\phi_{\sigma_*C'}^*}
\\
E_2^{h\aut(C')} \ar[r]_{y_{C'}^* \circ (\Frob_p^{-1})^*} &
E_2^{h\aut(\sigma_*C')}
}
$$
\end{proof}

\subsection{A resolution of $E_2^{h\Gamma^1_{Gal}}$}

In this section we explain how the statement of Theorem~\ref{thm:Q(ell)} 
changes when we replace the spectrum $E_2^{h\Gamma_{Gal}}$ with
the spectrum $E_2^{h\Gamma^1_{Gal}}$, where $\Gamma^1_{Gal}$ is the norm
$1$ subgroup defined in Section~\ref{sec:Galois}.

\begin{thm}\label{thm:hGamma1}
The spectrum $E_2^{h\Gamma^1_{Gal}}$ is equivalent to 
the homotopy fiber of the map
$$
\TMF_{K(2)} \times \TMF_{K(2)} \xrightarrow{p_2 \circ \phi_q^* - p_1 \circ
\phi_f^*} \TMF_0(\ell)_{K(2)} $$
where the maps $p_i$ are projections and the maps $\phi_q^*$ and $\phi_f^*$
are the maps defined in Section~\ref{sec:Q(ell)const}.
\end{thm}

The proof of Theorem~\ref{thm:hGamma1} follows the same lines as the proof
of Theorem~\ref{thm:Q(ell)}.  Namely, one uses Equation~\ref{eq:Jgamma} to
deduce the analog of Proposition~\ref{prop:Ehgammacosimp}: the spectrum
$E_2^{h\Gamma^1}$ is equivalent to the totalization of a
semi-cosimplicial spectrum of the form
$$ \prod_{C' \in X^{ss}} E_2^{h\aut(C')} \times E_2^{h\br{\aut(C')}}
\rightrightarrows \prod_{(C',H) \in X^{ss}_0(\ell)} E_2^{h\aut(C',H)}.
$$
One then forms the analog of Diagram~\ref{diag:topbottom}
\begin{equation}\label{diag:topbottom1}
\xymatrix@R+1em{
\prod_{C'} E_{C'}^{h\aut(C')} \times E_{C'}^{h\aut(C')}
\ar[d]|{(\phi_{C'}^*) \times ((\ell^{-r}\phi_{C'}\phi)^*)} 
\ar@<.5ex>[r] \ar@<-.5ex>[r] 
&
\prod_{(C',H)} E_{C'}^{h\aut(C',H)}
\ar[d]^{(\phi_{C'}^*)}
\\
\prod_{C'} E_2^{h\aut(C')} \times E_2^{h\br{\aut(C')}}
\ar@<.5ex>[r] \ar@<-.5ex>[r] 
&
\prod_{(C',H)} E_2^{h\aut(C',H)}
}
\end{equation}
where the coface maps of the top row correspond to $p_2 \circ \phi_q^*$ and 
$p_1 \circ \phi_f^*$, and the coface maps of the bottom row are determined
by Proposition~\ref{prop:Jface}.
Here $\phi$ is the endomorphism of $C$ of degree $\ell^{2r+1}$ that we
chose in Section~\ref{sec:orbits1}.
The essential point to the commutativity of Diagram~\ref{diag:topbottom1} is
the analog of Diagram~\ref{diag:essential}: for each $(C',H)$, 
the following diagram commutes.
$$
\xymatrix@C+1em{
E_{C_H}
\ar[r]^{\phi_H^*} 
\ar[d]_{(\ell^{-r}\phi_{C_H}\phi)^*}
&
E_{C'} \ar[d]^{\phi_{C'}^*}
\\
E_2
\ar[r]_{L_{g_{(C',H)}^1}}
&
E_2
}
$$
We deduce that there is an equivalence between $E_2^{h\Gamma^1}$ and the
homotopy fiber of the map
$$
\TMF_{K(2),\br{\FF}_p} \times \TMF_{K(2),\br{\FF}_p} 
\xrightarrow{p_2 \circ \phi_q^* - p_1 \circ
\phi_f^*} \TMF_0(\ell)_{K(2),\br{\FF}_p}.
$$
The Galois equivariance of this equivalence follows the same line of
verification that appears in the proof of Lemma~\ref{lem:verticalGal}, and
thus Theorem~\ref{thm:hGamma1} is obtained by taking Galois homotopy fixed
points.

\end{document}